\documentclass[11pt,reqno]{amsart}      
\usepackage[english]{babel}
\usepackage[utf8x]{inputenc}
\usepackage[T1]{fontenc}
\usepackage[letterpaper, margin=1.2in]{geometry}
\usepackage{mathtools}
\usepackage[dvipsnames]{xcolor}
\usepackage{hyperref}
\usepackage{dsfont}
\usepackage{enumerate}
\usepackage{dsfont}
\usepackage[pdftex]{graphicx}
\usepackage{caption}
\usepackage{subcaption}
\usepackage{tikz}
\usetikzlibrary{positioning,decorations.pathreplacing,shapes.misc}
\tikzset{cross/.style={cross out, draw=black, fill=none, minimum size=2*(#1-\pgflinewidth), inner sep=0pt, outer sep=0pt}, cross/.default={2pt}}

\allowdisplaybreaks

\usepackage{amsmath, amsthm, amssymb}
\usepackage{latexsym,graphicx,multicol,mathrsfs}

\newtheorem{theorem}{Theorem}[section]
\newtheorem{definition}[theorem]{Definition}
\newtheorem{proposition}[theorem]{Proposition}
\newtheorem{lemma}[theorem]{Lemma}
\newtheorem{corollary}[theorem]{Corollary}

\newtheorem{remark}[theorem]{Remark}

\newtheorem*{convention}{Convention}
\newtheorem*{proposition*}{Proposition}

\DeclareMathOperator{\var}{var}
\DeclareMathOperator{\cov}{cov}

\DeclareMathOperator{\Hess}{Hess}
\DeclareMathOperator{\Id}{Id}
\DeclareMathOperator{\supp}{supp}

\usepackage{autonum}

\begin{document}

\title[Hydrodynamic limit of the Kawasaki dynamics]{Hydrodynamic limit of the Kawasaki dynamics on the 1d-lattice with strong, finite-range interaction}

\author{Younghak Kwon}
\address{Department of Mathematics, University of California, Los Angeles}
\email{yhkwon@math.ucla.edu}

\author{Georg Menz}
\address{Department of Mathematics, University of California, Los Angeles}
\email{gmenz@math.ucla.edu}

\author{Kyeongsik Nam}
\address{Department of Mathematics, University of California, Los Angeles}
\email{ksnam@math.ucla.edu}


\date{\today}

\begin{abstract}

We derive the hydrodynamic limit of the Kawasaki dynamics for  the one-dimensional conservative system of unbounded real-valued spins with arbitrary strong, quadratic and finite-range interactions. This extends prior results for non-interacting spin systems. The result is obtained by adapting two scale approach of Grunewald, Otto, Villani and Westdickenberg combined with the authors' recent approach on  conservative systems with strong interactions.
\end{abstract}

\maketitle

\section{Introduction} \label{intro}

The broader scope of this article is the study of the continuum approximations of large discrete systems. Since the fundamental observation of Boltzmann that large particle systems in equilibrium are governed by Gibbs states, understanding the connection between discrete systems and their approximation to the continuum has been one of the main challenges in statistical physics. One of the most actively studied problem in this field is the hydrodynamic limit, which can be thought as a dynamical version of law of large numbers. This means that in a proper time and space macroscopic scales, the random evolution of microscopic system can be macroscopically described by a solution of deterministic partial differential equations. \\


In the 1980s, hydrodynamic limits were deduced for different settings. Two examples are the simple exclusion process (cf.~\cite{KOV89}) and the Kawasaki dynamics (cf.~\cite{Fri87}). In an attempt to provide more general strategy for deducing hydrodynamic limits, Guo, Papanicolaou and Varadhan proposed the martingale method to derive the hydrodynamic limit of the Kawasaki dynamics (cf.~\cite{GPV88}). In~\cite{Yau91}, Yau introduced the entropy method based on Gronwall-type estimate for a relative entropy functional to deduce the hydrodynamic limit of Kawasaki dynamics. One advantage of entropy method against martingale method is that, it is simpler and gives stronger results. However, it has a limitation that it assumes a stronger assumptions on the initial data. Another problem for both methods are not being very general yet to be applied to various settings. For more details, we refer to the classical reference on hydrodynamic limits~\cite{KipLan99}. \\

In this article, we study the hydrodynamic limit of Kawasaki dynamics for a one-dimensional lattice system of unbounded real-valued spins with arbitrary strong, finite-range interactions. Kawasaki dynamics is a stochastic dynamics preserving mean spins. Even if the case of absent interactions was solved in 1980s (e.g.~\cite{Fri87},~\cite{GPV88} and~\cite{Yau91}), up to authors' knowledge, there was no progress of problem until this work. To understand the obstacles and difficulties one encounters when studying the hydrodynamic limit, let us introduce the notion of grand- canonical and canonical ensemble. The grand canonical ensemble~$\mu_N$ is a probability measure on~$\mathbb{R}^N$ given by
\begin{align}
\mu_N(dx) : = \frac{1}{Z} \exp\left( - H(x) \right) dx.
\end{align} 
Here,~$Z$ denotes a generic normalization constant and~$H$ is the Hamiltonian of the system. Let us consider the~$N-1$ dimensional hyperplane~$X_{N,m}$ given by
\begin{align}
X_{N, m} : = \left \{ x \in \mathbb{R}^{N} : \ \frac{1}{N} \sum_{i=1}^{N} x_i =m  \right \} \subset \mathbb{R}^{N}.
\end{align}
The canonical ensemble~$\mu_{N,m}$ is the restriction of~$\mu_N$ to the~$N-1$ dimensional hyperplane~$X_{N,m}$:
\begin{align}
    \mu_{N,m}  (dx) : = \mu_N  \left(dx  \mid \frac{1}{N}\sum_{i =1}^{N} x_i =m \right) = \frac{1}{Z} \mathds{1}_{ \left\{ \frac{1}{N} \sum_{i =1}^{N} x_i =m \right\}}\left(x\right) \exp\left( - H(x) \right) \mathcal{L}^{N-1}(dx),
\end{align}
where~$\mathcal{L}^{N-1}(dx)$ denotes the~$N-1$ dimensional Hausdorff measure restricted to the hyperplane~$X_{N,m}$. In particular, the canonical ensemble is a stationary distribution of the Kawasaki dynamics. \\

Recalling that the hydrodynamic limit can be understood as a dynamical version of law of large numbers, it is obvious that the problem becomes relatively much easier if the underlying stationary distribution is a product measure; because this leads to independence of random variables. If the Hamiltonian~$H$ is non-interacting, we observe that the grand canonical ensemble is a product measure whereas the canonical ensemble is not due to the restriction to a hyperplane. This makes deducing hydrodynamic limit for the Kawasaki dynamics non-trivial even in the non-interactive Hamiltonian case. However, given the equivalence of ensembles (cf.~\cite{KLM19}) meaning that the canonical ensemble is equivalent to properly modified grand canonical ensemble, it is not surprising that one is able to deduce the hydrodynamic limit of Kawasaki dynamics in the case of a non-interactive Hamiltonian. \\

The problem becomes a lot more subtle if we consider interactions between spins within the Hamiltonian~$H$. In this case, even the grand canonical ensemble is not a product measure, making the task of deducing hydrodynamic limit for the Kawasaki dynamics even more challenging. Another difficulty of studying Kawasaki dynamics is that the lack of understanding of properties of the canonical ensemble with interacting Hamiltonian. For example, it is very important to show that on the microscopic scale there is a fast equilibration. This is closely connected to the logarithmic Sobolev inequality uniform in the system size, which is a sufficient condition for the fast equilibration. Other ingredients that can be useful are decay of correlations, strict convexity of the coarse-grained Hamiltonian and its convergence. \\

The deadlock was broken recently in~\cite{Me11}. There, the uniform LSI, decay of correlations, and strict convexity of coarse-grained Hamiltonian were solved in the case of weak interactive Hamiltonian. This provided an important tool for studying the hydrodynamic limit of a weakly interactive system, though the hydrodynamic limit was not deduced in this work. The case of weak interactions does not face the problems we faced in the study of arbitrary strong interactions, because one would expect that everything is close to the case of absent interaction. Indeed, the results were obtained by a perturbation argument, proving that weakly interactive system is close to a perturbed non-interactive system. This article would have prepared the ground to derive the hydrodynamic limit in the case of weak interactions. But instead, the authors chose to tackle the much harder problem of studying arbitrary strong interactions. \\

The situation is much harder in strongly interactive case and thus it required a series of articles to provide the same tools and insights about the structure of the canonical ensemble. This is a culmination of a series of works done by the authors. Recently, the authors provided better understandings of the canonical ensembles with strong finite-range interactions in a series of articles (see~\cite{KwMe18},~\cite{KwMe19a},~\cite{KwMe19b} and~\cite{KLM19}). These include, but not restricted to, equivalence of ensembles, decay of correlations and uniform LSI for the canonical ensemble. Those preparatory studies made it feasible to attack the problem of deducing the hydrodynamic limit in the case of strong interactions.  \\

While in principle it might be possible to adapt the martingale method (cf.~\cite{GPV88}) or entropy method (cf.~\cite{Yau91}), we chose the two-scale approach introduced by Grunewald, Otto, Villani and Westdickenberg (cf.~\cite{GrOtViWe09}). In~\cite{GrOtViWe09}, a general strategy for proving the hydrodynamic limit was derived via two-scale approach, applying Gronwall-type estimate and uniform logarithmic Sobolev inequality. We follow two-scale approach because it is a quantitative method whereas entropy and martingale methods are more qualitative in nature. A recent progress of~\cite{Jar18} enables the relative entropy method quantitative for a different process. As a consequence of our past work, proving quantitative hydrodynamic limit of Kawasaki dynamics via two scale approach became possible. However, in this article, we do not prove the quantitative hydrodynamic limit of Kawasaki dynamics as this approach will result in sub-optimal scaling of convergence and unnecessarily complicate our argument. Nevertheless, the quantitative hydrodynamic limit would be an important ingredient when studying fluctuations not starting in equilibrium (see e.g.~\cite{Jar18}). \\

One possible way improve the scaling of the the convergence would be to adapt two-scale approach with a more carefully chosen mesoscopic dynamics, as was done in~\cite{DMOW18} for the non-interactive case. The main difference to~\cite{GrOtViWe09} is that~\cite{DMOW18} introduces a mesoscopic dynamics as the Galerkin approximation of the macroscopic dynamics, while~\cite{GrOtViWe09} uses a projection onto piece-wise constant functions to define the mesoscopic scale. This approach using Galerkin approximation has an advantage of gaining regularity of the mesoscopic scale, resulting an optimal error estimate. It would be a challenging and interesting problem to extend this approach to the case of strongly interactive Hamiltonian. \\

Let us mention main challenges when applying the two-scale approach in the case of strong interactions. First of all, the convergence of the one-dimensional coarse-grained Hamiltonian should be handled. In case of non-interacting spin system, the local Cram\`er theorem implies that the one-dimensional coarse-grained Hamiltonian converges to the Cram\`er transform of a single-site potential. However, this is not true anymore under existence of strong interactions. Second, a uniform LSI should be extended from one block to multi blocks. That is, we consider the ensemble with conservation laws in each block and deduce the uniform LSI independent of block size, number of blocks, and the whole system size. Last, due to the strong finite-range interactions, the neighboring blocks are not independent anymore, resulting that the coarse-grained Hamiltonian is not a sum of one-dimensional coarse-grained Hamiltonians.  \\

To overcome the first difficulty, we recall that the local Cram\`er theorem implies the uniform convergence of one-dimensional coarse-grained Hamiltonian to a Legendre transform of the free energy of the grand canonical ensemble in the absence of interactions. Motivated by this, we first prove the convergence of the free energy of the grand canonical ensemble under the presence of strong interactions. In fact, we show that the sequence of (non-normalized) free energy  is sub-additive up to moment bounds. The moments are then compared with Gaussian moments, resulting bounds uniform on system size and depend only on the mean spin~$m$. Then we argue that the coarse-grained Hamiltonian converges to the Legendre transform of the limit of the free energy of the grand canonical ensemble. We refer to Section~\ref{section two scale} for more details. \\

The second difficulty, when deducing the multi-block LSI, is handled by applying a combination of the two-scale approach (cf.~\cite{GrOtViWe09}) and the Zegarlinski decomposition (cf.~\cite{Zeg96}). We decompose the lattice into two types of blocks~$\Lambda_1$ and~$\Lambda_2$ motivated by Zegarlinski's decomposition (cf. Figure~\ref{f_two_scale_decomposition_of_lattice}). Then the measure is decomposed into a conditional distribution conditioned on~$\Lambda_2$ and marginal distribution. By a careful choice of~$\Lambda_1$ and~$\Lambda_2$, the conditional distribution factorizes and thus the uniform LSI for the conditional distributions follows from a uniform LSI for the canonical ensemble  (cf.~\cite{KwMe19b}) and the Tensorization Principle. For the marginal distribution, we apply Otto-Reznikoff Criterion (see~\cite{OttRez07}) where interactions between blocks are controlled via decay of correlations. Then a usual two-scale argument for LSI combines the LSIs for conditional and marginal distributions and the uniform LSI for the original measure is obtained. For more details, we refer to Section~\ref{s_proof_generalized_lsi}. \\

For the last difficulty, we artificially introduce an auxiliary Hamiltonian~$H_{\text{aux}}$ where we remove the interactions between neighboring blocks. Removing the interactions makes each block independent, and as a consequence, the corresponding coarse-grained Hamiltonian of~$M$ blocks is decomposed into a sum of~$M$ coarse-grained Hamiltonians of single blocks. Because we assume finite range interactions, the number of interactions we remove is relatively small compared to the whole system size. Therefore, as expected, we prove that difference between the coarse-grained Hamiltonians arising from the formal Hamiltonian~$H$ and an auxiliary Hamiltonian~$H_{\text{aux}}$ goes to~$0$ as we increase the block size~$K$. This is well explained in Section~\ref{s_proof_hydro}. \\

Let us comment on open questions and problems:
\begin{itemize}

\item Instead for finite-range interaction, could one deduce similar results for infinite-range, algebraically decaying interactions? More precisely, is it possible to extend the results of~\cite{MeNi14} from the gce to the ce? If yes, is the same order of algebraic decay sufficient, i.e.~of the order~$2+ \varepsilon$, or does one need a higher order of decay? For solving this problem one would have to overcome several difficulties. For example, generalizing the equivalence of ensembles (see~\cite{KwMe18}) would need new work. Also, because we use ideas of the Zegarlinski method, the arguments of this article are restricted to the one-dimensional lattice with finite-range interaction. Applying our method to infinite-range interaction would yield a cyclic dependence of the different parameters. A possible alternative approach to this problem is to generalize the approach of~\cite{OttRez07,Me13,MeNi14} from the canonical ensemble to the grand canonical ensemble.\\[-1ex]

\item Is it possible to consider more general Hamiltonians? For example, our argument is based on the fact that the single-site potentials are perturbed quadratic, especially when we use the results of~\cite{KwMe18}. One would like to have general super-quadratic potentials as was for example used in~\cite{MeOt13}. \\[-1ex]

\item Is it possible to generalize the results to vector-valued spin systems? \\ [-1ex]
\end{itemize}

We conclude this Section by giving an overview over the article. In Section~\ref{s_setting_and_main_results} we introduce precise setting and present main results. In Section~\ref{section two scale}, we state key ingredients and prove several auxiliary results. In Section~\ref{s_proof_generalized_lsi} and Section~\ref{s_proof_strict_convexity_cg_hamiltonian}, two main ingredients uniform LSI and strict convexity of the coarse-grained Hamiltonian are proved, respectively. In Section~\ref{s_proof_hydro}, we give the proof of the main result of this article, namely, hydrodynamic limit of Kawasaki dynamics.

\section*{Conventions and Notation}

\begin{itemize}
\item The symbol~$T_{(k)}$ denotes the term that is given by the line~$(k)$.
\item We denote with~$0<C<\infty$ a generic uniform constant. This means that the actual value of~$C$ might change from line to line or even within a line.
\item Uniform means that a statement holds uniformly in the system size~$N$, the mean spin~$m$ and the external field~$s$.
\item $a \lesssim b$ denotes that there is a uniform constant~$C$ such that~$a \leq C b$.
\item $a \sim b$ means that~$a \lesssim b$ and~$b \lesssim a$.
\item $\mathcal{L}^{k}$ denotes the $k$-dimensional Hausdorff measure. If there is no cause of confusion we write~$\mathcal{L}$.
\item $Z$ is a generic normalization constant. It denotes the partition function of a measure.  
\item For each~$N \in \mathbb{N}$,~$[N]$ denotes the set~$\left\{ 1, \ldots N \right\}$.
\item For a vector~$x \in \mathbb{R}^{N}$ and a set~$A \subset [N]$,~$x^A \in \mathbb{R}^{A}$ denotes the vector $ (x^A)_{i} = x_i$ for all~$i \in A$.
\item For a vector~$x \in \mathbb{R}^{N}$ and a set~$A \subset [N]$,~$\bar{x}^{A} = x^{[N] \setminus A} \in \mathbb{R}^{[N] \setminus A}$ denotes the vector $(\bar{x}^A )_{i} = x_i$ for all~$i \in [N] \setminus A$
\item For a function~$f : \mathbb{R}^{N} \to \mathbb{C}$, we denote with~$\supp f = \{i_1, \cdots, i_k \}$ the minimal subset of~$[N]$ such that~$f(x) = f(x_{i_1}, \cdots, x_{i_k})$. 
\end{itemize}

\section{Setting and main results} \label{s_setting_and_main_results}

\subsection{The model}
The Gibbs measure we consider throughout the paper is a canonical ensemble with  strong interactions. The simplest case of canonical ensembles, where all of the interactions are removed, is considered in  ~\cite{GrOtViWe09} and~\cite{DMOW18}. The interactions we consider is strong in the sense that interactions are  beyond the perturbative regime.

Let us describe the precise model. Let~$\Lambda$ be the sublattice given by~$\Lambda = [N] = \{1, \cdots, N\}$. We consider a system of unbounded continuous spins on~$\Lambda$. The formal Hamiltonian $H = H_N :\mathbb{R}^{N} \to \mathbb{R}$ of the system is defined as
\begin{align}\label{e_d_hamiltonian}
H(x) = \sum_{i =1 }^{N} \left( \psi (x_i) + \frac{1}{2}\sum_{j : \ 1 \leq |j-i| \leq R } M_{ij}x_i x_j \right),
\end{align}
where~$\psi (z ) = \frac{1}{2} z^2 + \psi_b (z)$. For each~$i \in [N]$, we define~$M_{ii} : =1$ and set~$x_j =0$ for all~$j \notin [N]$. We also make the following assumptions:
\begin{itemize}
\item The function~$\psi_b: \mathbb{R} \to \mathbb{R}$ satisfies
 \begin{align}\label{e_nonconvexity_bounds_on_perturbation}
 |\psi_b|_{\infty} + |\psi'_b|_{\infty}  + |\psi''_b|_{\infty} < \infty. 
 \end{align}
It is best to imagine~$\psi$ as a double-well potential with quadratic growth at infinity (see Figure~\ref{f_double_well}).  
\begin{figure}[t]
\centering
\begin{tikzpicture}
      \draw[->] (-2.5,0) -- (2.5,0) node[right] {$x$};
      \draw[->] (0,-1) -- (0,2.8) node[above] {};
      \draw[scale=1,domain=-1.9:1.9,smooth,variable=\x,blue] plot ({\x},{.3*\x*\x*\x*\x-.7*\x*\x-.3*\x+.2});
      
\node[align=center, below, blue] at (1.9, 1.9) {$\psi(x)$};      
\end{tikzpicture}
\caption{Example of a single-site potential~$\psi$}\label{f_double_well}
\end{figure}
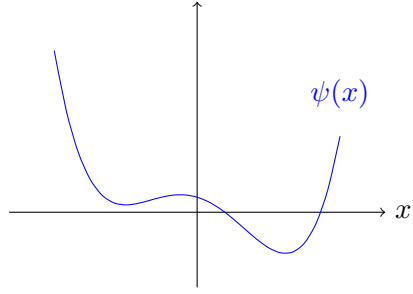

\item The interactions are symmetric, i.e., for any distinct~$i, j \in [N]$ with $|i-j| \leq R$,
\begin{align}
M_{ij} = M_{ji}. 
\end{align}

\item The fixed, finite number~$R \in \mathbb{N}$ models the range of interactions between the particles in the system i.e.~it holds that~$M_{ij}=0$ for all~$i,j$ such that~$|i-j|> R$. \\

\item The matrix~$(M_{ij})$ is strictly diagonal dominant i.e. for some~$\delta>0$, it holds for any~$i \in [N]$ that
\begin{align} \label{e_strictly_diagonal_dominant}
\sum_{ 1 \leq |j-i| \leq R} |M_{ij}| + \delta \leq M_{ii} = 1.
\end{align}

\item We assume spacial homogeneity of interactions. That is, there exists a function~$h : \mathbb{Z} \to (-1, 1)$ such that
\begin{align}
    M_{ij} = h(|i-j|) \qquad \text{for all distinct } i, j \in \mathbb{N}.
\end{align}
\end{itemize}

Let us define~$ X= X_{N,m}$ to be the~$(N-1)$-dimensional hyperplane with mean~$m$. More precisely, define
\begin{align}\label{e_d_X_lambda_m}
X = X_{N, m} : = \left \{ x \in \mathbb{R}^{N} : \ \frac{1}{N} \sum_{i=1}^{N} x_i =m  \right \} \subset \mathbb{R}^{N}.
\end{align}
We equip~$l^2$ inner product on~$X$ as follows:
\begin{align}
\langle x, \tilde{x} \rangle_{X} : = \sum_{i=1}^N x_i \tilde{x}_i.
\end{align}

The grand-canonical ensemble (gce)~$\mu_N$ associated to the Hamiltonian~$H$ is the probability measure on~$\mathbb{R}^{N}$ given by the Lebesgue density
\begin{align} \label{d_gc_ensemble}
\mu_N \left(dx\right) : = \frac{1}{Z} \exp\left( - H(x) \right) dx. 
\end{align}

The canonical ensemble (ce) emerges from the gce by conditioning on the mean spin
\begin{align} \label{e_spin_restriction}
  \frac{1}{N} \sum_{i =1}^{N} x_i = m.
\end{align}
More precisely, the ce~$\mu_{N,m} $ is the probability measure on~$X$
with density
\begin{align} 
\mu_{N,m}  (dx) : &= \mu_N  \left(dx  \mid \frac{1}{N}\sum_{i =1}^{N} x_i =m \right) \\
&= \frac{1}{Z} \mathds{1}_{ \left\{ \frac{1}{N} \sum_{i =1}^{N} x_i =m \right\}}\left(x\right) \exp\left( - H(x) \right) \mathcal{L}^{N-1}(dx), \label{d_ce}
\end{align}
where~$\mathcal{L}^{N-1}(dx)$ denotes the~$(N-1)$-dimensional Hausdorff measure supported on~$X$. \\

\subsection{Hydrodynamic limit of the Kawasaki dynamics} \label{s_main_hydro}

A natural dynamics for the conservative system is the Kawasaki dynamics, which is defined as follows.
Let~$A$ denote the second-order difference operator given by the~$N \times N$ matrix
\begin{align}
A_{ij} = N^2 \left(- \delta_{i, j-1} + 2 \delta_{i,j} - \delta_{i, j+1}\right),
\end{align}
where we define~$\delta_{i,0} = \delta_{i,N}$ and~$\delta_{i, N+1} = \delta_{i,1}$. The Kawasaki dynamics is a stochastic process~$X(t) \in \mathbb{R}^{N}$ satisfying the following stochastic differential equation:
\begin{align}
dX(t) = - A \nabla H(X(t)) dt + \sqrt{2A}dB(t),
\end{align}
where~$B(t)$ denotes a standard Brownian motion on~$\mathbb{R}^N$. The Kawasaki dynamic preserves its mean spins, i.e.,
\begin{align}
\frac{1}{N} \sum_{i=1}^N X_i (t) = \frac{1}{N} \sum_{i=1}^N X_i (0) = m.
\end{align}
This implies that we can restrict the state space~$\mathbb{R}^N$ to the hyperplane~$X= X_{N,m}$ and consider the corresponding ce~$\mu_{N,m}$. If the process~$X_t$ is distributed according to~$f \mu_{N,m}$, then the time dependent probability density~$f = f(t,x)$ satisfies
\begin{align} \label{e_kawasaki_de}
\frac{\partial}{\partial t} (f \mu_{N, m}) = \nabla \cdot \left( A \nabla f \mu_{N,m} \right).
\end{align}

In order to define a continuous counterpart of the configuration space $X_{N,m}$,
let us define the space~$\bar{X}$ of piecewise constant, mean~$m$ functions on~$\mathbb{T}^1 = \mathbb{R} \backslash \mathbb{Z}$ by
\begin{align} \label{e_def_barX}
\bar{X} : = \left\{ \bar{x} : \mathbb{T}^1 \to \mathbb{R} ; \ \bar{x} \text{ is constant on } \left( \frac{j-1}{N}, \frac{j}{N} \right] \text{ for } j=1, \cdots, N, \text{ and has mean } m \right\}.
\end{align}
We shall identify the space~$X= X_{N,m}$ with~$\bar{X}$ by the following relation:
\begin{itemize}
\item For each~$x \in X$, the step function~$\bar{x} \in \bar{X}$ associated to~$x$ is
\begin{align} \label{e_X_to_barX_identification}
\bar{x}(\theta) = x_j, \qquad \text{if } \theta \in \left( \frac{j-1}{N}, \frac{j}{N} \right].
\end{align}
\item For each step function~$\bar{x} \in \bar{X}$, the corresponding vector~$x \in X$ is
\begin{align} \label{e_barX_to_X_identification}
x_j = \bar{x} \left( \frac{j}{N} \right), \qquad j=1, \cdots, N.
\end{align}
\end{itemize}

\medskip

We equip the space of locally integrable functions~$f : \mathbb{T}^1 \to \mathbb{R}$ and has mean~$m$ with~$H^{-1}$ norm by
\begin{align}
\| f \|_{H^{-1}}^2 = \int_{\mathbb{T}^1} \omega^2 (\theta) d\theta,
\end{align}
where~$\omega$ is a function such that
\begin{align}
\omega ' = f, \qquad \int \omega(\theta) d\theta = 0.
\end{align}

\medskip

Now we are ready to formulate our main result, namely the hydrodynamic limit of the Kawasaki dynamics. We establish that the evolution along the Kawasaki dynamics  gets close to the solution to  a certain nonlinear parabolic equation as $N\rightarrow \infty$.

\begin{theorem}\label{p_hydrodynamic_limit} Let~$f= f(t,x)$ be a solution of the Kawasaki dynamics~\eqref{e_kawasaki_de} with initial condition~$f(0, \cdot) = f_0 (\cdot)$. Assume that there is a positive constant~$C>0$ such that 
\begin{align} \label{e_entropy_f0}
\int f_0 (x) \log f_0 (x) \mu_{N,m} (dx) \leq CN. 
\end{align}
Assume also that there is a~$\zeta_0 \in L^2 (\mathbb{T}^1)$ such that~$\int \zeta_0 d\theta = m  $ and
\begin{align} \label{e_barx_zeta0}
\lim_{N \to \infty} \int \| \bar{x} - \zeta_0 \|_{H^{-1}}^2 f_0 (x) \mu_{N,m} (dx) =0.
\end{align}

Let~$\zeta = \zeta(t, \theta)$ be the unique weak solution of the nonlinear parabolic equation
\begin{align} \label{e_heat_eqn}
\begin{cases} \frac{\partial \zeta}{\partial t} = \frac{\partial^2}{\partial \theta^2} \varphi ' (\zeta), \\ \zeta(0, \cdot) = \zeta_0,
\end{cases}
\end{align}
where $\varphi$ is defined as 
\begin{align}
    \varphi(m) := \lim_{N\rightarrow \infty}   -\frac{1}{N} \log  \int_{ \{\frac{1}{N}\sum_{i=1}^{N} x_i = m \} } \exp\left(-H(x)\right) \mathcal{L}^{N-1}(dx) \label{211}
\end{align}
  Then for any~$T>0$, it holds that
\begin{align}
\lim_{N \to \infty} \sup_{0 \leq t \leq T} \int \| \bar{x} - \zeta (t, \cdot) \|_{H^{-1}}^2 f(t,x) \mu_{N,m} (dx) =0.
\end{align}
\end{theorem}

Here, we say that  $\zeta = \zeta(t, \theta)$ is a weak solution of~\eqref{e_heat_eqn} on~$[0, T] \times \mathbb{T}^1$ if
\begin{align}
\zeta \in L_t^{\infty}(L_{\theta}^2), \qquad \frac{\partial \zeta}{\partial t} \in L_t ^2 (H_{\theta}^{-1}), \qquad \text{and} \qquad  \varphi' (\zeta) \in L_t ^2 (L_{\theta}^2),
\end{align}
and
\begin{align}
\left \langle \xi, \frac{\partial \zeta}{\partial t} \right \rangle_{H^{-1}} = - \int_{\mathbb{T}^1} \xi \varphi' (\zeta) d\theta \qquad \xi \in L^2, \text{for almost every }  t\in [0,T]  
\end{align}

The proof of Theorem~\ref{p_hydrodynamic_limit} is given in Section~\ref{s_proof_hydro}.

\begin{remark} 
There are some issues in Theorem \ref{p_hydrodynamic_limit} to be resolved. First, one has to verify that the pointwise limit of \eqref{211} exists and   is differentiable. This will be  established in Section \ref{section two scale}. In addition, the existence and uniqueness of a weak solution of~\eqref{e_heat_eqn} follows from the standard argument in the nonlinear PDE theory (see for example~\cite[Lemma 38]{GrOtViWe09}).
\end{remark}

\begin{remark} \label{remark 2.3}
The quantity inside the limit of \eqref{211}, denoted by $\bar{H}_N(m)$, represents the distribution $f_N(m)dm$  of the mean value $(x_1+\cdots+x_N)/N$  under $\mu_{N,m}$:
\begin{align}
    f_N(m)dm = \frac{1}{Z_N}e^{-N \bar{H}_N(m) }dm.
\end{align} 
In the case of ce without interactions, i.e. $M_{ij}=0$, as a consequence of local Cram\`er theorem (see \cite[Proposition 31]{GrOtViWe09}),  $\varphi$ in \eqref{211} is a Legendre transform of the logarithmic generating function of the distribution $\frac{1}{Z}e^{-\psi(x)}dx$. On the other hand, in the presence of interactions, $\varphi$ in \eqref{211} can also be expressed in terms of the
Legendre transform of the thermodynamic free energy. This point will be discussed   in  Section \ref{section two scale}.
\end{remark}



\section{Two-scale decomposition} \label{section two scale}
In this section, we introduce a two-scale decomposition method, originally introduced in \cite{GrOtViWe09},  which plays a   crucial  role to study the concentration properties of the ce and their hydrodynamic limit. Then, we state key results on the logarithmic Sobolev inequality and the strict convexity for the coarse-grained Hamiltonian, generalizing the previous results in~\cite{KwMe19b},~\cite{KwMe18}, which are crucial to implement a two-scale approach to establish a hydrodynamic limit.

\medskip

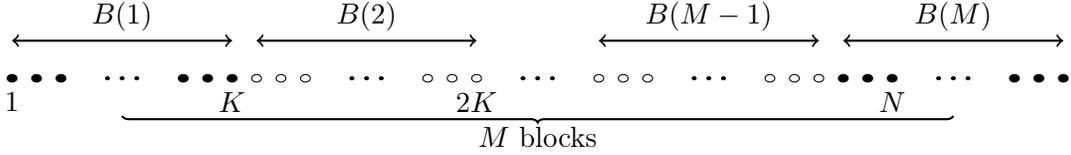
\begin{figure}[t]
\centering
\begin{tikzpicture}[xscale=1.3]

\draw[fill] (0,0) circle [radius=0.05];
\draw[fill] (.25,0) circle [radius=0.05];
\draw[fill] (.5,0) circle [radius=0.05];

\draw[fill] (0.975,0) circle [radius=0.02];
\draw[fill] (1.125,0) circle [radius=0.02];
\draw[fill] (1.275,0) circle [radius=0.02];

\draw[fill] (1.75,0) circle [radius=0.05];
\draw[fill] (2,0) circle [radius=0.05];
\draw[fill] (2.25,0) circle [radius=0.05];

\draw (2.5, 0) circle [radius=0.05];
\draw (2.75,0) circle [radius=0.05];
\draw (3,0) circle [radius=0.05];

\draw[fill] (3.475,0) circle [radius=0.02];
\draw[fill] (3.625,0) circle [radius=0.02];
\draw[fill] (3.775,0) circle [radius=0.02];

\draw (4.25,0) circle [radius=0.05];
\draw (4.5,0) circle [radius=0.05];
\draw (4.75,0) circle [radius=0.05];

\draw[fill] (5.225,0) circle [radius=0.02];
\draw[fill] (5.375,0) circle [radius=0.02];
\draw[fill] (5.525,0) circle [radius=0.02];

\draw (6,0) circle [radius=0.05];
\draw (6.25,0) circle [radius=0.05];
\draw (6.5,0) circle [radius=0.05];

\draw[fill] (6.975,0) circle [radius=0.02];
\draw[fill] (7.125,0) circle [radius=0.02];
\draw[fill] (7.275,0) circle [radius=0.02];

\draw (7.75,0) circle [radius=0.05];
\draw (8,0) circle [radius=0.05];
\draw (8.25,0) circle [radius=0.05];

\draw[fill] (8.5, 0) circle [radius=0.05];
\draw[fill] (8.75, 0) circle [radius=0.05];
\draw[fill] (9, 0) circle [radius=0.05];

\draw[fill] (9.475,0) circle [radius=0.02];
\draw[fill] (9.625,0) circle [radius=0.02];
\draw[fill] (9.775,0) circle [radius=0.02];

\draw[fill] (10.25, 0) circle [radius=0.05];
\draw[fill] (10.5, 0) circle [radius=0.05];
\draw[fill] (10.75, 0) circle [radius=0.05];

\node[align=center, below] at (0,-.05) {$1$};
\node[align=center, below] at (2.25,-.05) {$K$};
\node[align=center, below] at (4.75,-.05) {$2K$};
\node[align=center, below] at (9,-.05) {$N$};
\node[align=center, below] at (4.5,-1.2) {};

\draw[decorate,decoration={brace,mirror},thick] (1.125,-.5) -- node[below]{$M$ blocks} (9.625,-.5);

\draw[thick,<->] (0,.5) -- (2.25,.5);
\node[align=center, above] at (1.125,.5) {$B(1)$};
\draw[thick,<->] (2.5, .5) -- (4.75,.5);
\node[align=center, above] at (3.625, .5) {$B(2)$};

\draw[thick,<->] (6, .5) -- (8.25,.5);
\node[align=center, above] at (7.125, .5) {$B(M-1)$};

\draw[thick,<->] (8.5, .5) -- (10.75,.5);
\node[align=center, above] at (9.625, .5) {$B(M)$};

\end{tikzpicture}
\caption{Block decomposition of lattice~$[N]$.}\label{f_blcok_decomposition}
\end{figure}

Let us divide~$N$ spins into~$M$ blocks with size~$K$ (see Figure~\ref{f_blcok_decomposition}), denoted by
\begin{align}
    B(l) : = \{ (l-1) K +1, \cdots, lK \} \qquad \text{for each } l \in \{ 1, \cdots, M \}.
\end{align}
We then define the mesoscopic space~$Y$ as
\begin{align} \label{d_mesoscopic_y}
Y= Y_{M, m} = \left \{  (y_1, \cdots, y_M ) ; \  \frac{1}{M} \sum_{l=1}^M y_l = m  \right \}.
\end{align}

The~$L^2$ inner product in~$Y$ is defined as follows:
\begin{align}
\langle y, \tilde{y} \rangle_Y = \frac{1}{M} \sum_{l=1}^M y_l \tilde{y}_l.
\end{align}

The projection $P = P_{N,K} : X \to Y$ is defined via
\begin{align}
P (x_1, \cdots, x_N ) = \left( y_1, \cdots, y_M  \right), \qquad y_l = \frac{1}{K} \sum_{i  \in B(l)} x_i.
\end{align}
We observe that the adjoint operator~$P^* : Y \to X$ given by
\begin{align}
P^* (y_1, \cdots, y_M) = \frac{1}{N} ( \underbrace{ y_1, \cdots, y_1}_{K \text{ times}}, \cdots, \underbrace{y_M, \cdots, y_M}_{K \text{ times}} ) 
\end{align}
satisfies the identity~$PN P^* = \Id_Y$, where~$\Id_Y$ is the identity operator on~$Y$.

\begin{remark}
For notational simplicity, we assumed that all blocks~$B(l)$ have equal size~$K$ so that~$N=MK$. If~$N/K$ is not an integer, we decompose~$[N]$ into~$M$ blocks with different sizes~$K_1, \cdots, K_M$. More precisely, we define
\begin{align}
    Px = (y_1, \cdots, y_M ),
\end{align}
where
\begin{align}
y_l = \frac{1}{K_l} \sum_{i \in B(l)} x_i \qquad \text{for each } l \in \{ 1, \cdots, M \}.
\end{align}
The space~$Y$ is defined as
\begin{align}
Y = \left \{  (y_1, \cdots, y_M ) ; \  \frac{1}{M} \sum_{l=1}^M \alpha_l y_l = m, \text{ where } \alpha_l = \frac{MK_l}{N} \right \}.
\end{align}
Here, we choose block sizes~$\{K_l\}_{l=1}^{M}$ carefully so that~$1 \leq \alpha_l \leq 2$ for all~$l \in \{1, \cdots, M\}$. 
\end{remark}

Finally, we  disintegrate the ce~$\mu_{N,m}$ into the conditional measure~$\mu_{N,m}(dx | y ) = \mu_{N,m}(dx | Px=y)$ and the marginal measure~$\bar{\mu}_{N,m}(y)$ defined on $Y$. This means that for any test function~$\xi$,
\begin{align}  \label{e_decomposition_of_measures}
\int \xi d\mu_{N,m} = \int_Y \left( \int \xi(x) \mu_{N,m}(dx | y ) \right) \bar{\mu}_{N,m} (dy).
\end{align}

\subsection{Key ingredients: logarithmic Sobolev inequality  and coarse-grained Hamiltonian}

The two-scale decomposition method has been successfully used to study the hydrodynamic limit of the Kawasaki dynamics of ce without interactions.  Key ingredients to establish hydrodynamic limit are the uniform  logarithmic Sobolev inequality for the conditional distributions and the strict convexity of the coarse-grained Hamiltonian (see Section  \ref{s_proof_hydro} for details).
In this section, we state new results on the logarithmic Sobolev inequality and coarse-grained Hamiltonians in the context of ce with strong interactions.

\medskip

Let us first introduce the definition of the logarithmic Sobolev inequality (LSI):
\begin{definition}[Logarithmic Sobolev Inequality (LSI)]
Let~$X$ be a Euclidean space. A Borel probability measure~$\mu$ satisfies a logarithmic Sobolev inequality with constant~$\varrho>0$ if for all test functions~$f \geq 0$,
\begin{align}
\int f \log f d\mu  - \int f d \mu \log \left( \int f d\mu \right) \leq \frac{1}{2\varrho} \int \frac{|\nabla f|^2}{f} d\mu.
\end{align}
When~$X = \mathbb{R}^N$, we say~$\mu$ satisfies a uniform LSI with constant~$\varrho>0$ if~$\varrho$ is independent of the system size~$N$. 
\end{definition}

There have been numerous works on studying LSI for the conservative spin systems. Important works include  \cite{LuYa93}, where  a martingale method was implemented, and \cite{GrOtViWe09}, where a two-scale method was introduced. Recently, the uniform LSI for~$\mu_{N,m}$, the ce with strong interactions, was obtained 
in~\cite{KwMe19b}:

\begin{lemma}[Theorem 2 in~\cite{KwMe19b}] \label{p_uniform_lsi} 
The ce~$\mu_{N,m}$ given by~\eqref{d_ce} satisfies a uniform LSI($\varrho$), where~$\varrho>0$ is independent of the system size~$N$, the external field~$s$ and the mean spin~$m \in \mathbb{R}$.
\end{lemma}

Let us recall that the measure~$\mu_{N,m}$ conditions on the mean value~$m$ of the spins~$x_1, \ldots, x_N$. We therefore call the LSI for the measure~$\mu_{N,m}$ the \emph{one-block LSI}. \\

The uniform LSI implies a strong concentration property of the corresponding Gibbs measure and  provides an exponential decay of entropy along the associated dynamics. More precisely, the uniform LSI for the~$\mu_{N,m}$ implies the exponential decay in the relative entropy along the Kawasaki dynamics: 
\begin{align} \label{e_rel_entropy_lsi}
\text{Ent}(f(t, x) \mu_{N,m} | \mu_{N,m} ) \leq \exp ( -CN^{-2}t) \text{Ent}(f(0, x)\mu_{N,m}  | \mu_{N,m} ),
\end{align}
where~$\text{Ent}$ denotes the relative entropy function.  We refer to~\cite[Remark 3]{Me11} for more details. \\

The first main result of this section is an improvement of Lemma \ref{p_uniform_lsi} to multi blocks. More precisely, let us recall that the measure~$\mu_{N,m}(dx|y)$ conditions on the mean spins~$y_1, \ldots y_M$ on each block. We show that this measure also satisfies the uniform LSI, which is called \emph{multi-block LSI}. This is one of the key ingredients to implement the two-scale approach to establish the hydrodynamic limit. It is also is highly non-trivial because, due to the interactions between blocks, the measure~$\mu_{N,m}(dx|y)$ does not tensorize.  

\begin{theorem}\label{p_generalized_lsi} The conditional measure~$\mu_{N,m} (dx | y)$ satisfies a uniform LSI($\varrho$), where~$\varrho>0$ is independent of the system size~$N$, the external field~$s$, the mean spin~$m$, and the macroscopic state~$y$.
\end{theorem}

Theorem \ref{p_generalized_lsi} will be proved in Section \ref{s_proof_generalized_lsi}.\\

Now, we define and study some properties about the coarse-grained Hamiltonian. 
Recall the disintegration
\begin{align}
 \mu_{N,m}(dy) =  \mu_{N,m}(dx|y)  \bar{\mu}_{N,m}(dy).
\end{align}
The coarse-grained Hamiltonian $\bar{H}_Y(y)$ is defined to be a Hamiltonian corresponding to   $\bar{\mu}_{N,m}(dy)$:
\begin{align}
    \bar{\mu}_{N,m}(dy) = \exp(-N \bar{H}_Y(y))dy.
\end{align}
In other words, one can define a coarse-grained Hamiltonian~$\bar{H}_Y : Y \to \mathbb{R}$ as follows:
\begin{align}
\bar{H}_Y(y) &: = - \frac{1}{N} \log \int_{Px = y} \exp(-H(x)) \mathcal{L}^{N-M}(dx) \label{e_def_cg_hamiltonian}.
\end{align}
In particular, in the simple case $Px = (1/N) \sum_{i=1}^N x_i$, we define~$\bar{H}_N : \mathbb{R} \to \mathbb{R}$  by
\begin{align}
    \bar{H}_N(m) &: = -\frac{1}{N} \log  \int_{ \{\frac{1}{N}\sum_{i=1}^{N} x_i = m \} } \exp\left(-H(x)\right) \mathcal{L}^{N-1}(dx) \label{e_def_1d_cg_ham}.
\end{align}

The strict convexity of coarse-grained Hamiltonian of the ce plays a crucial role in~\cite{GrOtViWe09} to establish a uniform LSI and the hydrodynamic limit without interactions. In the one-block setting, this has been verified for strong interactions $\mu_{N,m}$ (see~\cite{KwMe18} or~\cite{KwMe19b}).

\begin{lemma}[Lemma 1 in~\cite{KwMe19b}] \label{l_strict_convex_1dim_cg}
The coarse-grained Hamiltonian~$\bar{H}_N : X \to \mathbb{R}$ is uniformly strictly convex. In other words, there exists a   constant~$C$ such that for any $N\geq 1$ and $m\in \mathbb{R}$,
\begin{align}
\frac{1}{C} \leq \bar{H}_N '' (m) \leq C .
\end{align}
\end{lemma}

In the second main result of this section we extend this result to the multi-block case:

\begin{theorem} \label{p_strict_convexity_cg_hamiltonian} The coarse-grained Hamiltonian~$\bar{H}_Y$ is uniformly strictly convex. In other words, there exists   a constant~$\lambda>0$ independent of the system size~$N$, the external field~$s$, and the mean spin~$m$, such that for any $y\in Y$,
\begin{align}
\lambda \Id_{Y} \leq \Hess_Y \bar{H}_Y(y) \leq  \frac{1}{\lambda} \Id_{Y} .
\end{align}
\end{theorem}

Theorem \ref{p_strict_convexity_cg_hamiltonian} will be proved in Section \ref{s_proof_strict_convexity_cg_hamiltonian}.

\begin{remark} One should compare Theorem~\ref{p_generalized_lsi} and Theorem~\ref{p_strict_convexity_cg_hamiltonian} with Lemma~\ref{p_uniform_lsi} and Lemma~\ref{l_strict_convex_1dim_cg}, respectively. In Lemma~\ref{p_uniform_lsi} and Lemma~\ref{l_strict_convex_1dim_cg}, the authors considered the case where Gibbs measure has only one constraint
\begin{align}
\frac{1}{N} \sum_{i=1}^{N} x_i = m.
\end{align}
The setting of Theorem~\ref{p_generalized_lsi} and Theorem~\ref{p_strict_convexity_cg_hamiltonian} is more general in the sense that the measure~$\mu_{N,m}(dx|y)$ has multiple constraints, having one conservation law for each block:
\begin{align}
\frac{1}{K} \sum_{i \in B(l)} x_i = y_l, \qquad l =1, \cdots, M.
\end{align}
That is, if we let~$M=1$, the statements of Theorem~\ref{p_generalized_lsi} and Theorem~\ref{p_strict_convexity_cg_hamiltonian} reduce to that of Lemma~\ref{p_uniform_lsi} and Lemma~\ref{l_strict_convex_1dim_cg}, respectively.
\end{remark}

Finally, we verify that the pointwise limit of the one-dimensional coarse grained Hamiltonian $\bar{H}_N(m)$ in  \eqref{e_def_1d_cg_ham} exists as the system size goes to infinity. It turns out that this limit, denoted by $\varphi$, is a function that appears in the nonlinear parabolic equation \eqref{e_heat_eqn}. The following lemma provides a quantitative convergence of $\bar{H}_N(m)$ as $N\rightarrow \infty$.

\begin{proposition} \label{l_convergence_cg_ham}
There exists a differentiable function~$\varphi : \mathbb{R} \to \mathbb{R}$ such that for each $m\in \mathbb{R}$,
\begin{align} 
    \bar{H}_N(m) \to \varphi(m) \qquad \text{as } N \to \infty.  \label{e_convergence_cg_ham_limit}
\end{align}
Moreover, there exist a positive constant~$C$   such that for any $N\geq 1$ and $m\in \mathbb{R}$,
\begin{align}
\left| \bar{H}_N (m) - \varphi (m)  \right| \leq C \frac{m^2 +1}{N}.
\end{align}
\end{proposition}

 \begin{remark}
 In the case of ce without interactions, as mentioned in Remark \ref{remark 2.3}, $\varphi$ is a Legendre transform of the logarithmic generating function of the distribution $\frac{1}{Z}e^{-\psi(x)}dx$, and moreover the convergence in  \eqref{e_convergence_cg_ham_limit}
 holds in the $C^2$ topology.  This is a consequence of the local Cram\`er theorem obtained in \cite[Proposition 31]{GrOtViWe09}. However, Cram\`er's large deviation principle for the mean $(x_1+\cdots+x_N)/N$ does not hold in general under the presence of dependencies among random variables. In our case, i.e. in the presence of strong interactions, we find a candidate $\varphi$ for Theorem \ref{p_hydrodynamic_limit} as a limit of one-dimensional coarse-grained Hamiltonians. It turns out that  $\varphi$ can also be represented by a Legendre transform of the thermodynamic free energy of the corresponding gce.
 \end{remark}

For the rest of section, we prove Proposition \ref{l_convergence_cg_ham}. As a key ingredient, we first establish a sharp moment estimates with respect to the ce~$\mu_{N,m}$.
Then, equipped with moment estimates, we establish Proposition \ref{l_convergence_cg_ham} by showing that $\varphi$ is  a Legendre transform of the thermodynamic free energy of gce.

\subsection{Moment estimates} \label{section moment}
In this section, we obtain sharp moment estimate for the ce and gce. We first study moments for the gce, and then analyze in the case of ce using the principle of equivalence of ensembles.

\subsubsection{Moment estimates under grand canonical ensembles}

We establish a sharp moment estimate under the following gce with external fields
\begin{align} \label{d_modified_gce}
\mu^{\sigma}_N (dx) : = \frac{1}{Z} \exp\left( \sigma \sum_{i=1}^N x_i - H(x) \right) dx.
\end{align}
Next lemma provides the first moment bound under the gce \eqref{d_modified_gce}.
\begin{lemma} \label{l_gce_moment_estimate}
For any $N\geq 1$ and $i\in [N]$, we have
\begin{align} \label{e_gce_moment_estimate}
\left|    \mathbb{E}_{\mu_{N}^{\sigma}} \left[ X_i \right]  \right| \lesssim \sigma +1 .
\end{align}
\end{lemma}
\medskip

It is delicate to estimate the first moment directly under the measure \eqref{d_modified_gce}. We overcome this problem by comparing the first moment under the gce and Gaussian ensemble. Since it is straightforward to compute the first moment under the Gaussian measure, one can finally deduce Lemma \ref{l_gce_moment_estimate}.
\\

\noindent \emph{Proof of Lemma~\ref{l_gce_moment_estimate}.} \ As mentioned above, proof consists of the following two steps:
\begin{itemize}
\item Transfer from the gce to Gaussian ensembles using interpolation.
\item Obtain a sharp estimate on the first moment under the Gaussian ensembles.
\end{itemize}
\medskip

\textbf{Step 1.} Comparison with Gaussian ensembles. \\

For $s \in\mathbb{R}$, let us define Hamiltonian
\begin{align}
H_s(x) :=  \sum_{i=1}^N \left( \frac{1}{2}x_i^2 +  \sum_{j:1\leq |j-i|\leq R} M_{ij}x_ix_j  + s \psi_b (x_i) \right)
\end{align}
and the corresponding measure
\begin{align}
\nu_{N,s}^{\sigma}(dx): = \frac{1}{Z} \exp\left(\sigma \sum_{i=1}^{N} x_i -H_s(x)\right)dx.
\end{align}
In particular, we note that
\begin{align}
    \nu_{N, 1} ^{\sigma} = \mu_N ^{\sigma}
\end{align}
and~$\nu_{N, 0}^{\sigma}$ is a Gaussian ensemble. Let us fix~$i \in [N]$. We now apply an interpolation to obtain
\begin{align}
    \mathbb{E}_{\nu_{N, 1} ^{\sigma}} \left[ X_i \right] - \mathbb{E}_{\nu_{N, 0}^{\sigma}} \left[ X_i \right] &= \int_0 ^1 \frac{d}{ds} \mathbb{E}_{\nu_{n, s}^{\sigma}} \left[ X_i \right] ds \\
    & = \int_0 ^1 \cov_{\nu_{N,s}^{\sigma}}\left( X_i, -\sum_{j=1}^{N} \psi_b (x_j) \right) ds \\
    & = - \int_0 ^1 \sum_{j=1}^{N} \cov_{\nu_{N, s}^{\sigma}} \left( X_i, \psi_b (x_j ) \right) ds.
\end{align}
Because~$\nu_{N, s}$ is again a gce, we have by Theorem~\ref{p_decay_of_correlations_gce} that
\begin{align}
    \left| \cov_{\nu_{N, s}^{\sigma}} \left(X_i, \psi_b (x_j) \right) \right| &\leq C_s \| \nabla X_i \|_{L^2 (\nu_{N, s}^{\sigma})} \| \nabla \psi_b (x_j) \|_{L^2 (\nu_{N, s}^{\sigma})} \exp \left( -C_s |i-j| \right) \\
    & \leq C_s \| \psi_b ' \|_{\infty} \exp\left(-C_s |i-j| \right).
\end{align}
Because the constant~$C_s$ uniformly bounded by a constant~$C>0$ for~$s \in [0, 1]$, we have
\begin{align}
    \left| \mathbb{E}_{\nu_{N, 1} ^{\sigma}} \left[ X_i \right] - \mathbb{E}_{\nu_{N, 0}^{\sigma}} \left[ X_i \right] \right| \leq C \sum_{j=1}^{N} \exp \left( -C |i-j|\right) \lesssim 1.
\end{align}
That is,
\begin{align} \label{e_bound_by_gaussian}
    \left| \mathbb{E}_{\mu_N ^{\sigma}} \left[ X_i \right] \right| = \left| \mathbb{E}_{\nu_{N, 1} ^{\sigma}} \left[ X_i \right] \right| \lesssim \left| \mathbb{E}_{\nu_{N, 0} ^{\sigma}} \left[ X_i \right] \right| +1
\end{align}

\medskip
\textbf{Step 2.} First moment estimate under Gaussian ensembles. \\

Note that for any $\kappa\in \mathbb{R}$,
the Gaussian measure $\nu_{N,0}^{\sigma}$ has a mean  vector $(\eta_1^\sigma,\cdots,\eta_N^\sigma)$ given by
\begin{align} 
(\eta_1^\sigma,\cdots,\eta_N^\sigma)^T =  M^{-1} ( \sigma ,\cdots, \sigma)^T.
\end{align}
Because the Hamiltonian~$H_s(x)$ is strictly convex, the measure~$\nu_{N, 0}^{\sigma}$ satisfies a uniform LSI independent of~$\sigma$ and thus we have the exponential decay of correlations (see~\cite{HeMe16} for example). For multivariate Gaussian distribution, the covariance matrix is given by the inverse of the quadratic coefficient matrix~$M$, i.e.~$\cov_{\nu_{N, 0}^{\sigma}} (X) = M^{-1} $. In particular, the coefficient of~$M^{-1}$ decays exponentially in the sense that
\begin{align}
    (M^{-1})_{ij} \leq C \exp \left(- C |i-j| \right) \qquad \text{for each } i, j.
\end{align}
Therefore there exist constants $c_1',c_2'>0$  such that for any $i=1,\cdots,N$,
\begin{align}\label{h12}
| \mathbb{E}_{ \nu_{N,0}^{\sigma}} x_i     | \leq  c_1'\sigma + c_2'.
\end{align} 

\medskip

Thus, by~\eqref{e_bound_by_gaussian} and~\eqref{h12},  one can deduce that  there exist constants $c_1,c_2  >0$  
such that for any   $\sigma\in \mathbb{R}$,
\begin{align}
\left|  \mathbb{E}_{ \mu_N^\sigma} \left[x_i\right] \right|  \leq c_1\sigma + c_2.
\end{align}

\qed

Because we have a sharp bound for the variances by Poincar\`e inequality, one can deduce the following corollary as a consequence of Lemma \ref{l_gce_moment_estimate}. \\

\begin{corollary} \label{l_gce_second_moment}
For any $N\geq 1$ and $i\in [N]$,
\begin{align}
\left|    \mathbb{E}_{\mu_{N}^{\sigma}} \left[ X_i ^2 \right]  \right| \lesssim \sigma^2 +1.
\end{align}
\end{corollary}

\medskip

  \noindent \emph{Proof of Lemma~\ref{l_gce_second_moment}}
  By Poincar\`e inequality,
  \begin{align}
     \var(X_i) \lesssim 1.
  \end{align}
 Combining this with Lemma \ref{l_gce_moment_estimate}, we conclude the proof.

\subsubsection{Moment estimates under canonical ensembles}
Using the moment estimates under gce and the principle of equivalence of observables (Proposition~\ref{p_equivalence_observables}), we obtain the moment estimate for the ce.

\begin{lemma} \label{l_ce_moment_estimate}
For any $N\geq 1$ and $i\in [N]$,
\begin{align}
\left|    \mathbb{E}_{\mu_{N,m}} \left[ X_i \right]  \right| \lesssim m +1 .
\end{align}
\end{lemma}

\medskip

\noindent \emph{Proof of Lemma~\ref{l_ce_moment_estimate}.} \ We first prove that there exist constants $\gamma_1,\gamma_2,\gamma'_1,\gamma'_2\in \mathbb{R}$ such that
for any $\sigma$ and $m$ satisfying
\begin{align}
\frac{d}{d\sigma}A(\sigma) = m,
\end{align}
we have
\begin{align}  \label{h1}
\gamma_1m+\gamma_2\leq  \sigma \leq \gamma_1'm+\gamma_2'.
\end{align}
By Lemma~\ref{l_strict_convex_free_energy}, there exists  $C>0$ such that for any $\sigma\in \mathbb{R}$,
\begin{align} \label{h2}
 \frac{1}{C} \leq \frac{d^2}{d\sigma^2}A(\sigma)  \leq C.
\end{align}
Also,  note that
\begin{align} \label{h3}
\frac{d}{d\sigma}A(0) = \frac{1}{N} \mathbb{E}_{\mu_N^{0}} \left[ \sum_{i=1}^N X_i \right].
\end{align}
This quantity is uniformly bounded in $N$ thanks to Lemma \ref{l_gce_moment_estimate}.
Thus, by \eqref{h2} and \eqref{h3}, we have \eqref{h1} for some constants $\gamma_1,\gamma_2,\gamma'_1,\gamma'_2\in \mathbb{R}$. \\

By the equivalence of observable result (Proposition~\ref{p_equivalence_observables}), we have
\begin{align} \label{h4}
\left| \mathbb{E}_{ \mu_{N,m}} \left[x_i\right] - \mathbb{E}_{ \mu_N^\sigma} \left[x_i\right] \right|  = O(\frac{1}{N}).
\end{align}
Thus, by Lemma \ref{l_gce_moment_estimate}, \eqref{h1},  and \eqref{h4}, proof is concluded. \\ 
\qed 

Because the ce~$\mu_{N,m}$ satisfies a uniform LSI and hence Poincar\'e inequality, the variance of the ce is well behaved. Therefore we have the following statement. \\

\begin{corollary} \label{l_ce_second_moment}
For any $N\geq 1$ and $i\in [N]$,
\begin{align}
\left|    \mathbb{E}_{\mu_{N, m}} \left[ X_i ^2 \right]  \right| \lesssim m^2 +1 .
\end{align}
\end{corollary}

\medskip

\subsection{Convergence of the coarse-grained Hamiltonian}

Let us define a (non-normalized) free energy of the gce
\begin{align}
    a_N(\sigma) : = \log \int_{\mathbb{R}^N}  \exp \left( \sigma \sum_{i=1}^{N} x_i -H_N(x) \right) dx.
\end{align}
First of all, we prove that for each~$\sigma \in \mathbb{R}$,~$a_N(\sigma)$ and~$a_N ' (\sigma)$ are sub-additive up to constants.

\begin{lemma} \label{l_subadditive_gce}
There exists a positive constant~$C$ such that for any~$N_1,N_2 \in \mathbb{N}$ and~$\sigma \in \mathbb{R}$,
\begin{align} 
\left| a_{N_1+N_2}(\sigma) - a_{N_1}(\sigma) - a_{N_2} (\sigma)\right| &\leq C(\sigma^2 +1), \label{e_subadditive_gce}\\
\left| a_{N_1+N_2} ' (\sigma) - a_{N_1} ' (\sigma) - a_{N_2} ' (\sigma)\right| &\leq C(|\sigma| +1). \label{e_subadditive_derivative_gce}
\end{align}
\end{lemma}

\noindent \emph{Proof of~\eqref{e_subadditive_gce} in Lemma~\ref{l_subadditive_gce}.} \ We write
\begin{align}
&a_{N_1+N_2}(\sigma) - a_{N_1}(\sigma) - a_{N_2}  (\sigma) \\
& \qquad = \log \frac{ \int_{\mathbb{R}^{N_1+N_2}}  \exp \left( \sigma \sum_{k=1}^{N_1+N_2} w_k -H_{N_1+N_2}(w) \right) dw}{  \int_{\mathbb{R}^{N_1}}  \exp \left( \sigma \sum_{i=1}^{N_1} u_i -H_{N_1}(u) \right) du \cdot \int_{\mathbb{R}^{N_2}}  \exp \left( \sigma \sum_{j=1}^{N_2} v_i -H_{N_2}(v) \right) dv   } \\
& \qquad = \log \frac{ \int_{\mathbb{R}^{N_1+N_2}}  \exp \left(\sigma \sum_{k=1}^{N_1+N_2} w_k -H_{N_1+N_2}(w) \right) dw}{  \int_{\mathbb{R}^{N_1+N_2}}  \exp \left( \sigma \sum_{i=1}^{N_1} u_i + \sigma \sum_{j=1}^{N_2} v_j -H_{N_1}(u)  -H_{N_2}(v) \right) dudv   }.
\end{align}
Let us denote
\begin{align} \label{d_inm}
I_{N_1,N_2} = \{ (i, j) \ | \  i \in \{1, \cdots, N_1\}, j \in \{N_1+1, \cdots, N_1+N_2 \}, |i-j| \leq R \}.
\end{align}

Writing~$w= (u, v) \in \mathbb{R}^N \times \mathbb{R}^M$, we have
\begin{align} \label{e_diff_block_hamiltonian}
H_{N_1+N_2}(w) - H_{N_1}(u) - H_{N_2} (v) =   \sum_{(i,j) \in I_{N_1,N_2}}  M_{ij} u_i v_j.
\end{align}
Thus we can write
\begin{align}
    a_{N_1+N_2}(\sigma) - a_{N_1}(\sigma) - a_{N_2}(\sigma)  =\log \left( \mathbb{E}_{\mu_{N_1} ^{\sigma}  \bigotimes \mu_{N_2} ^{\sigma}} \left[ \exp \left( -  \sum_{(i, j) \in I_{N_1,N_2}}  M_{ij} u_i v_j \right) \right] \right). \label{e_subadditive}
\end{align}

We shall only prove that~\eqref{e_subadditive} is bounded from above as the proof of lower bound is almost identical to that of upper bound. \\

To begin with, an application of Young's inequality yields
\begin{align}
T_{\eqref{e_subadditive}}
&  \leq \log \left( \mathbb{E}_{\mu_{N_1} ^{\sigma} \bigotimes \mu_{N_2} ^{\sigma} } \left[ \exp \left( \frac{1}{2}  \sum_{(i, j) \in I_{N_1,N_2}}  |M_{ij}| \left(u_i^2 + v_j ^2 \right) \right) \right]\right) \\
& = \log \left( \mathbb{E}_{\mu_{N_1} ^{\sigma}} \left[ \exp \left( \frac{1}{2}  \sum_{(i, j) \in I_{N_1,N_2}}  |M_{ij}| u_i^2  \right) \right] \right) \label{e_subadditive_upper1} \\
&  \quad + \log \left( 
\mathbb{E}_{\mu_{N_2} ^{\sigma}} \left[ \exp \left( \frac{1}{2}  \sum_{(i, j) \in I_{N_1,N_2}}  |M_{ij}| v_j^2  \right) \right] \right).  \label{e_subadditive_upper2}
\end{align}

We then apply an interpolation and get
\begin{align}
T_{\eqref{e_subadditive_upper1}} &= \log\left( \int_{\mathbb{R}^{N_1}} \exp \left( \sigma \sum_{i=1}^{N_1} u_i -H_{N_1}(u) + \frac{1}{2}  \sum_{(i, j) \in I_{N_1,N_2}}  |M_{ij}| u_i^2  \right) du \right) \\
& \quad - \log \left( \int_{\mathbb{R}^N_1} \exp \left( \sigma \sum_{i=1}^{N_1} u_i - H_{N_1}(u) \right) du \right) \\
& = \int_0 ^1 \frac{d}{ds} \log\left( \int_{\mathbb{R}^N_1} \exp \left( \sigma \sum_{i=1}^{N_1} u_i -H_{N_1}(u) + s \cdot \frac{1}{2}  \sum_{(i, j) \in I_{N_1, N_2}}  |M_{ij}| u_i^2  \right) du \right) ds \\
& = \int_0^1 \mathbb{E}_{\mu_{N_1} ^{\sigma}(s)} \left[  \frac{1}{2}  \sum_{(i, j) \in I_{N_1,N_2}}  |M_{ij}| u_i^2 \right] ds, \label{e_subadditive_upper1_equiv}
\end{align}
where~$\mu_{N_1} ^{\sigma} (s)$ is the probability distribution given by
\begin{align}
\mu_{N_1} ^{\sigma} (s) (dx) := \frac{1}{Z} \exp \left( \sigma \sum_{i=1}^{N_1} x_i  -H_{N_1}(x) + s \cdot \frac{1}{2}  \sum_{(i, j) \in I_{N_1,N_2}}  |M_{ij}| x_i^2 \right) dx. 
\end{align}
We observe that
\begin{align}
&H_{N_1}(x) -  s \cdot \frac{1}{2}  \sum_{(i, j) \in I_{N_1,N_2}}  |M_{ij}| x_i^2 
\\
& \qquad = \sum_{i =1 }^{N_1-R} \left( \psi (x_i) + s_i x_i +\frac{1}{2}\sum_{j : \ 1 \leq |j-i| \leq R } M_{ij}x_i x_j \right) \\
& \qquad \quad + \sum_{i=N_1-R+1}^{N_1} \left( \psi(x_i) - s \cdot \frac{1}{2} \sum_{\substack{j \in \{N_1+1, \cdots, N_1+N_2 \} \\ 1 \leq |j-i| \leq R}} |M_{ij}| x_i^2  + s_i x_i + \frac{1}{2}\sum_{\substack{j \in \{1, \cdots, N_1 \} \\ 1 \leq |j-i| \leq R}} M_{ij}x_i x_j \right)
\end{align}
Due to the strictly diagonal dominant assumption~\eqref{e_strictly_diagonal_dominant},
\begin{align}
\frac{1}{2} - s \cdot \frac{1}{2} \sum_{\substack{j \in \{N_1+1, \cdots, N_1+N_2 \} \\ 1 \leq |j-i| \leq R}} |M_{ij}|  \geq \frac{1}{2}\sum_{\substack{j \in \{1, \cdots, N_1 \} \\ 1 \leq |j-i| \leq R}} M_{ij} + \frac{1}{2} \delta.
\end{align}
This means that the interaction terms of~$\mu_N^{\sigma}(s)$ also satisfy the strictly diagonal dominant assumption~\eqref{e_strictly_diagonal_dominant} and hence~$\mu_N^{\sigma}(s)$ is a gce. Therefore an application of Corollary~\ref{l_gce_second_moment} implies
\begin{align}
    \left| T_{\eqref{e_subadditive_upper1}} \right| & \overset{\eqref{e_subadditive_upper1_equiv}}{\leq} \int_0 ^1  \frac{1}{2} \sum_{(i, j) \in I_{N_1,N_2}} |M_{ij}| \sup_{s \in [0, 1]} \left( \mathbb{E}_{\mu_{N_1} ^{\sigma}(s)} \left[ u_i ^2 \right] \right) ds \\
    & \overset{Corollary~\ref{l_gce_second_moment}}{\lesssim} \frac{1}{2} \sum_{(i, j) \in I_{N_1,N_2}} |M_{ij}| (\sigma^2 +1) \sim \sigma ^2 +1
\end{align}
Similarly, one gets~$|T_{\eqref{e_subadditive_upper2}}| \lesssim \sigma^2 +1$ and thus
\begin{align}
    T_{\eqref{e_subadditive}} \lesssim \sigma^2 +1.
\end{align}
\medskip

\noindent \emph{Proof of~\eqref{e_subadditive_derivative_gce} in Lemma~\ref{l_subadditive_gce}.} \ Following the notations used in the proof of~\eqref{e_subadditive_gce} in Lemma~\ref{l_subadditive_gce}, we write
\begin{align}
&a_{N_1+N_2} ' (\sigma) - a_{N_1} ' (\sigma) - a_{N_2} ' (\sigma) \\
& \qquad = \mathbb{E}_{\mu_{N_1+N_2}^{\sigma}} \left[ \sum_{k=1}^{N_1+N_2} w_k \right] - \mathbb{E}_{\mu_{N_1}^{\sigma}} \left[ \sum_{i=1}^{N_1} u_i \right] -\mathbb{E}_{\mu_{N_2}^{\sigma}} \left[ \sum_{j=1}^{N_2} v_j \right] \\
& \qquad = \mathbb{E}_{\mu_{N_1+N_2}^{\sigma}} \left[ \sum_{k=1}^{N_1+N_2} w_k \right] - \mathbb{E}_{\mu_{N_1}^{\sigma} \bigotimes \mu_{N_2}^{\sigma}} \left[ \sum_{k=1}^{N_1+N_2} w_k \right] \label{e_n+m_tensor_diff}
\end{align}
Let us recall the definition~\eqref{d_inm} of~$I_{N_1,N_2 }$. For~$s \in [0,1]$, define
\begin{align}
    \nu_{N_1+N_2}^{\sigma} (s) (dx) : = \frac{1}{Z} \exp \left( \sigma \sum_{k=1}^{N_1+N_2} w_k - H_{N_1+N_2}(w) + s \sum_{(i, j) \in I_{N_1, N_2}} M_{ij} u_i v_j \right).
\end{align}
In particular, recalling~\eqref{e_diff_block_hamiltonian}, it holds that
\begin{align}
    \nu_{N_1+N_2}^{\sigma} (0) = \mu_{N_1+N_2}^{\sigma}, \qquad \nu_{N_1+N_2}^{\sigma} (1)  = \mu_{N_1} ^{\sigma} \otimes \mu_{N_2} ^{\sigma}.
\end{align}
Therefore we apply interpolation to obtain
\begin{align}
    \left| T_{\eqref{e_n+m_tensor_diff}}\right|  & = \left| -  \int_0 ^1 \frac{d}{ds} \mathbb{E}_{\nu_{N_1+N_2 }^{\sigma}(s)}\left[ \sum_{k=1}^{N_1+N_2} w_k \right] dx \right| \\
    & \leq \int_0 ^1 \left| \cov_{\nu_{N_1+N_2}^{\sigma}(s)} \left( \sum_{k=1}^{N_1+N_2} w_k , \sum_{(i, j) \in I_{N_1,N_2 }} M_{ij} u_i v_j \right)  \right| ds \\
    & \leq \int_0 ^1 \sum_{k=1}^{N_1+N_2} \left|\cov_{\nu_{N_1+N_2}^{\sigma}(s)} \left(  w_k , \sum_{(i, j) \in I_{N_1,N_2}} M_{ij} u_i v_j \right)  \right| ds. \label{e_diff_interpolation_cov}
\end{align}
We observe that~$\nu_{N_1+N_2}^{\sigma}(s)$ is a gce and hence Theorem~\ref{p_decay_of_correlations_gce} holds. Because~$|I_{N_1, N_2}| \leq 2R^2$, an application of the decay of correlation (Theorem~\ref{p_decay_of_correlations_gce}) and second moment estimate (Lemma~\ref{l_gce_second_moment}) yields
\begin{align}
\left|\cov_{\nu_{N_1+N_2}^{\sigma}(s)} \left(  w_k , \sum_{(i, j) \in I_{N_1,N_2}} M_{ij} u_i v_j \right)  \right| &\leq C_s \left( 2R^2 (\sigma^2 +1) \right) ^{\frac{1}{2}} \exp \left( -C_s \text{dist}(k, I_{N_1,N_2}) \right)  \\
& \leq C (|\sigma| +1) \exp\left( -C \text{dist}(k, I_{N_1,N_2}) \right) \label{e_covariance_estimate}
\end{align} 
for some~$C>0$. Plugging~\eqref{e_covariance_estimate} into~\eqref{e_diff_interpolation_cov} yields
\begin{align}
\left| T_{\eqref{e_n+m_tensor_diff}}\right| \leq C(|\sigma|+1) \sum_{k=1}^{N_1+N_2} \exp\left( -C \text{dist}(k, I_{N_1,N_2}) \right)  = C\left( |\sigma| +1 \right).
\end{align}
\qed
\medskip

In particular, an application of Fekete's Lemma implies the following statement. \\

\begin{corollary}
The normalized free energy
\begin{align}
    A_{N}(\sigma) := \frac{1}{N} a_N = \frac{1}{N} \log \int_{\mathbb{R}^N}  \exp \left( \sigma \sum_{i=1}^{N} x_i - H_N (x) \right) dx
\end{align}
and its derivative~$A_N ' $ converge pointwise to some functions~$A, B :\mathbb{R} \to \mathbb{R}$ as~$N \to \infty$, respectively.
\end{corollary}
\medskip

Next, we provide quantitative bounds of the convergences of~$A_N$ and~$A_N '$. \\

\begin{lemma} \label{l_quantitative_conv_An}
There exists a positive constant~$C$   such that
\begin{align}
    |A_N (\sigma) - A (\sigma) | &\leq C \frac{\sigma^2 +1}{N} \qquad \text{for all } \sigma \in \mathbb{R}, \label{e_quantitative_conv_an} \\
     |A_N '(\sigma) - B(\sigma) | &\leq C \frac{|\sigma| +1}{N} \qquad \text{for all } \sigma \in \mathbb{R}. \label{e_quantitative_conv_an_der}
\end{align}
\end{lemma}
\medskip

\noindent \emph{Proof of Lemma~\ref{l_quantitative_conv_An}.} \ We shall only provide the proof of~\eqref{e_quantitative_conv_an} as there is only a cosmetic difference between the proof of~\eqref{e_quantitative_conv_an} and that of~\eqref{e_quantitative_conv_an_der}. \\

Let us fix~$N \in \mathbb{N}$. We claim that for each~$k \in \mathbb{N}$,
\begin{align} \label{e_induction_claim}
    \left| A_{kN} (\sigma) - A_N (\sigma) \right| \leq C \frac{k-1}{k} \cdot \frac{\sigma^2 +1}{N}.
\end{align}
First of all,~\eqref{e_induction_claim} is obviously true for~$k=1$. The case~$k=2$ also holds by putting~$M=N$ in~\eqref{e_subadditive_gce} and dividing it by~$2N$. Let us assume that~\eqref{e_induction_claim} holds for some~$p \in \mathbb{N}$. That is,
\begin{align} \label{e_induction_hypo}
    \left| A_{pN} (\sigma) - A_N (\sigma) \right| \leq C \frac{p-1}{p} \cdot \frac{\sigma^2 +1}{N}.
\end{align}
Then
\begin{align}
\left| A_{(p+1)N} (\sigma) - A_N (\sigma) \right| &= \left| \frac{a_{(p+1)N} (\sigma) }{(p+1)N} - \frac{a_N (\sigma) }{N}   \right| \\
& \leq \left| \frac{a_{(p+1)N}(\sigma) - a_{pN} (\sigma) - a_N (\sigma)}{(p+1)N}  \right| + \left| \frac{a_{pN}(\sigma) - p a_N (\sigma) }{(p+1)N}  \right| \\
& \overset{\eqref{e_subadditive_gce}}{\leq} C \frac{1}{p+1} \cdot \frac{\sigma^2 +1}{N} + \frac{p}{p+1} \left|A_{pN}(\sigma) - A_N (\sigma) \right|  \\
& \overset{\eqref{e_induction_hypo}}{\leq} C \frac{1}{p+1} \cdot \frac{\sigma^2 +1}{N}  + C \frac{p-1}{p+1} \cdot \frac{\sigma^2 +1}{N} = C \frac{p}{p+1} \cdot \frac{ \sigma^2 +1}{N}.
\end{align}
Therefore~\eqref{e_induction_claim} holds for~$k=p+1$ as well and thus it holds for all~$k \in \mathbb{N}$. We now take~$k \to \infty$ in~\eqref{e_induction_claim} to conclude that
\begin{align}
\left| A(\sigma) - A_N (\sigma) \right| \leq C \frac{\sigma^2 +1}{N}.
\end{align}
\qed

\medskip

Lemma~\ref{l_quantitative_conv_An} implies that if restricted to any closed interval [a, b], the convergence of~$A_N ' $ is uniform. Moreover,~$A_N'$ is differentiable for each~$N$ and in particular, is continuous. Thus we have the following statement: \\

\begin{corollary} \label{l_derivatives_match}
The function $A$ is a~$C^1$ function, and $A'=B$. In other words,
\begin{align}
    A' (\sigma) = \lim_{N \to \infty} A_N ' (\sigma) \qquad \text{for each } \sigma \in \mathbb{R}.
\end{align}
\end{corollary}
\medskip

Let~$\mathcal{H}_N$,~$\varphi : \mathbb{R} \to \mathbb{R}$ be the Legendre transforms of~$A_N$ and~$A$, respectively. That is,
\begin{align}
 \mathcal{H}_N (m) & := \sup_{\sigma \in \mathbb{R}} \left( \sigma m - A_N (\sigma) \right), \\
\varphi(m) & := \sup_{\sigma \in \mathbb{R}} \left( \sigma m - A(\sigma) \right). \label{d_varphi}
\end{align}
We recall that
\begin{align}
A_N ' (\sigma) = \frac{1}{N} \mathbb{E}_{\mu_N ^{\sigma}} \left[ \sum_{i=1}^{N} X_i \right]  \qquad \text{and} \qquad  A_N '' (\sigma) & = \frac{1}{N} \var_{\mu_N ^{\sigma}} \left( \sum_{i=1}^{N} X_i \right).
\end{align}
Then Lemma~\ref{l_gce_moment_estimate} and Lemma~\ref{l_variance_estimate} imply that there exists a uniform positive constant~$C$ such that
\begin{align} \label{e_an_derivatives_bound}
-C \left( 1 + |\sigma| \right) \leq A_N '(\sigma) \leq C \left( 1+ |\sigma| \right) \qquad \text{and} \qquad 
\frac{1}{C}  \leq A_N ''(\sigma) \leq C .
\end{align}
Because the bounds~\eqref{e_an_derivatives_bound} are uniform on~$N$, it also holds that
\begin{align} \label{e_a_derivative_bound}
-C \left( 1 + |\sigma| \right) \leq A'(\sigma) \leq C \left( 1+ |\sigma| \right)
\end{align}
and~$A$ is strictly convex in the sense that
\begin{align} \label{e_strict_convexity_A}
\frac{1}{C} (x-y) \leq A'(x) - A'(y) \leq C (x-y) \qquad \text{for any } x \geq y. 
\end{align}
The strict convexity of~$A_N$ implies that for each~$N$, there exists a unique real number~$\sigma_N \in \mathbb{R}$ such that
\begin{align} \label{d_sigman}
\mathcal{H}_N (m) = \sup_{\sigma \in \mathbb{R}} \left( \sigma m - A_N (\sigma) \right) = \sigma_N m - A_N (\sigma_N).
\end{align}
We also denote~$\sigma_{\infty}$ by the unique real number satisfying
\begin{align} \label{d_sigmainfty}
\varphi (m) = \sup_{\sigma \in \mathbb{R}} \left( \sigma m - A (\sigma) \right) = \sigma_{\infty} m - A (\sigma_{\infty}).
\end{align}

Next, we prove that~$\mathcal{H}_N \to \varphi$ pointwise as~$N \to \infty$. \\

\begin{lemma} \label{l_varphi_est}
There exists a positive constant~$C$  such that
\begin{align}
\left| \mathcal{H}_N (m) - \varphi(m) \right| \leq C \frac{m^2 +1}{N} \qquad \text{for all } m \in \mathbb{R}.
\end{align}
\end{lemma}
\medskip

\noindent \emph{Proof of Lemma~\ref{l_varphi_est}.} \ By definition~\eqref{d_sigman} and~\eqref{d_sigmainfty} of~$\sigma_N$ and~$\sigma_{\infty}$, it holds that
\begin{align} \label{e_sigman_sigmainfty_relation}
A_N ' (\sigma_N ) = A ' (\sigma_{\infty}) = m.
\end{align}
Let us recall the uniform linear bounds~\eqref{e_an_derivatives_bound},~\eqref{e_a_derivative_bound} of~$A_N '$,~$A'$ and the strict convexity~\eqref{e_an_derivatives_bound},~\eqref{e_strict_convexity_A} of~$A_N$,~$A$. They imply that there exist constants~$\gamma_1,\gamma_2,\gamma'_1,\gamma'_2\in \mathbb{R}$ with
\begin{align} \label{e_linear_bounds_sigma_m}
    \gamma_1 m + \gamma_2 \leq \sigma_N, \sigma_{\infty} \leq \gamma_1 ' m + \gamma_2 ' \qquad \text{for all } N \in \mathbb{N}.
\end{align}

Let us write
\begin{align}
\left|\mathcal{H}_N (m) - \varphi(m)\right| & \leq \left| \sigma_N - \sigma_{\infty} \right| |m| + \left| A_N (\sigma_N) - A(\sigma_{\infty}) \right|. \label{e_bound_hn_varphi}
\end{align}
Let us begin with the estimation of the first term in the right hand side of~\eqref{e_bound_hn_varphi}. Rearranging~\eqref{e_sigman_sigmainfty_relation} gives
\begin{align}
A' (\sigma_{\infty}) - A' (\sigma_N) = A_N ' (\sigma_N) - A' (\sigma_N)
\end{align}
Then the strict convexity~\eqref{e_strict_convexity_A} of~$A$ implies that there is a positive constant~$C$ with
\begin{align}
\frac{1}{C} |\sigma_{\infty} - \sigma_N | \leq \left| A' (\sigma_{\infty}) - A'(\sigma_N) \right| \leq C \left|\sigma_{\infty} - \sigma_N \right|.
\end{align}
Therefore we have
\begin{align} \label{e_bound_diff_sigman_sigmainfty}
\left| \sigma_{\infty} - \sigma_N \right| \leq C \left|A'(\sigma_{\infty}) - A'(\sigma_N ) \right| = C \left| A_N ' (\sigma_N ) - A' (\sigma_N ) \right| 
& \overset{\eqref{e_quantitative_conv_an_der}}{\leq} C \frac{(|\sigma_N|  +1)}{N} \overset{\eqref{h1}}{\leq} C \frac{ |m| +1}{N}. 
\end{align}
Let us turn to the estimation of the second term in the right hand side of~\eqref{e_bound_hn_varphi}. It holds that
\begin{align}
\left| A_N (\sigma_N) - A(\sigma_{\infty}) \right| & \leq \left| A_N (\sigma_N) - A(\sigma_N) \right| + \left| A(\sigma_N) - A(\sigma_{\infty}) \right| \\
& \overset{\eqref{e_quantitative_conv_an}}{\leq} C \frac{\sigma_N ^2 +1}{N} + \left|A(\sigma_N) - A(\sigma_{\infty}) \right| \\
& = C \frac{\sigma_N ^2 +1}{N} + \left| A'(\sigma_N^*) \right| \left| \sigma_N - \sigma_{\infty} \right|, \label{e_taylor}
\end{align}
where~$\sigma_N ^*$ is a real number between~$\sigma_N$ and~$\sigma_{\infty}$. Therefore a combination of~\eqref{e_linear_bounds_sigma_m},~\eqref{e_a_derivative_bound} and~\eqref{e_bound_diff_sigman_sigmainfty} yields
\begin{align} \label{e_taylor_res}
T_{\eqref{e_taylor}} \leq C \frac{m^2 +1}{N}.
\end{align}
Plugging the estimates~\eqref{e_bound_diff_sigman_sigmainfty} and~\eqref{e_taylor_res} into~\eqref{e_bound_hn_varphi} gives the desired estimate
\begin{align}
    \left| \mathcal{H}_N (m) - \varphi(m) \right| \leq C \frac{m^2 +1}{N}.
\end{align}
\qed

\medskip

The last ingredient for proving Lemma~\ref{l_convergence_cg_ham} is the local Cram\'er theorem. \\ 

\begin{lemma}[Theorem 2.6 in~\cite{KwMe18}] \label{l_cramer}
There exists a positive constant~$C$ such that for~$n$ large enough,
\begin{align}
 \left|   \bar{H}_N (m)  - \mathcal{H}_N (m) \right| \leq C \frac{1}{N} \qquad \text{for all } m \in \mathbb{R}.
\end{align}
\end{lemma}

We can now conclude the  proof of Proposition~\ref{l_convergence_cg_ham}. \\

\noindent \emph{Proof of Proposition~\ref{l_convergence_cg_ham}.} \ Let~$\varphi$ be the function defined by~\eqref{d_varphi}. A combination of Lemma~\ref{l_varphi_est} and Lemma~\ref{l_cramer} implies that
\begin{align}
    \left| \bar{H}_N(m) - \varphi(m) \right| & \leq \left| \bar{H}_N (m) - \mathcal{H}_N (m) \right| + \left| \mathcal{H}_N (m) - \varphi (m)\right| \\
    & \leq C \frac{1}{N} + C \frac{m^2 + 1}{N}  = C \frac{m^2 +1}{N}.
\end{align}
The differentiability of $\varphi$ is obvious. In fact, the Legendre transform of a~$C^1$ strictly convex function with superlinear growth is also differentiable.

\qed

\section{Logarithmic Sobolev inequality: proof of Theorem~\ref{p_generalized_lsi}} \label{s_proof_generalized_lsi}

The proof of Theorem~\ref{p_generalized_lsi} is motivated by Zegarlinski's decomposition which was used to prove the uniform LSI for the gce~$\mu_N$ (cf.~\cite{Zeg96}). In~\cite{KwMe19b}, the authors used this idea combined with the two-scale approach (cf.~\cite{GrOtViWe09}) to prove that the ce~$\mu_{N,m}$ satisfies a uniform LSI on one-dimensional lattice. In this proof, we adapt this idea using Zegarlinski's decomposition and two-scale approach to deduce the uniform LSI for the measure~$\mu_{N,m} (dx|y)$. \\

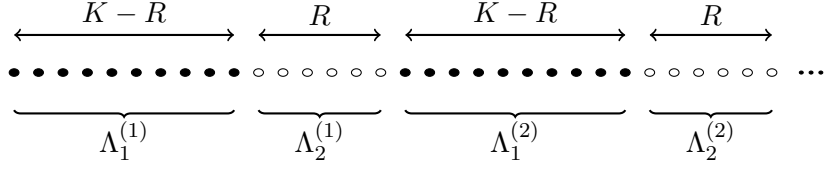
\begin{figure}[t]
\centering
\begin{tikzpicture}[xscale=1.3]

\draw[fill] (.5,0) circle [radius=0.05];
\draw[fill] (.75,0) circle [radius=0.05];
\draw[fill] (1,0) circle [radius=0.05];
\draw[fill] (1.25,0) circle [radius=0.05];
\draw[fill] (1.5,0) circle [radius=0.05];
\draw[fill] (1.75,0) circle [radius=0.05];
\draw[fill] (2,0) circle[radius=0.05];
\draw[fill] (2.25,0) circle[radius=0.05];
\draw[fill] (2.5,0) circle[radius=0.05];
\draw[fill] (2.75,0) circle[radius=0.05];

\draw (3,0) circle[radius=0.05];
\draw (3.25,0) circle[radius=0.05];
\draw (3.5,0) circle[radius=0.05];
\draw (3.75,0) circle[radius=0.05];
\draw (4,0) circle[radius=0.05];
\draw (4.25,0) circle[radius=0.05];

\draw[fill] (4.5,0) circle [radius=0.05];
\draw[fill] (4.75,0) circle [radius=0.05];
\draw[fill] (5,0) circle [radius=0.05];
\draw[fill] (5.25,0) circle [radius=0.05];
\draw[fill] (5.5,0) circle [radius=0.05];
\draw[fill] (5.75,0) circle [radius=0.05];
\draw[fill] (6,0) circle[radius=0.05];
\draw[fill] (6.25,0) circle[radius=0.05];
\draw[fill] (6.5,0) circle[radius=0.05];
\draw[fill] (6.75,0) circle[radius=0.05];

\draw (7,0) circle[radius=0.05];
\draw (7.25,0) circle[radius=0.05];
\draw (7.5,0) circle[radius=0.05];
\draw (7.75,0) circle[radius=0.05];
\draw (8,0) circle[radius=0.05];
\draw (8.25,0) circle[radius=0.05];

\draw[fill] (8.55,0) circle [radius=0.02];
\draw[fill] (8.65,0) circle [radius=0.02];
\draw[fill] (8.75,0) circle [radius=0.02];

\draw[decorate,decoration={brace,mirror},thick] (0.5,-.5) -- node[below]{$\Lambda_1 ^{(1)}$} (2.75,-.5);
\draw[decorate,decoration={brace,mirror},thick] (3,-.5) -- node[below]{$\Lambda_2 ^{(1)}$} (4.25,-.5);
\draw[decorate,decoration={brace,mirror},thick] (4.5,-.5) -- node[below]{$\Lambda_1 ^{(2)}$} (6.75,-.5);
\draw[decorate,decoration={brace,mirror},thick] (7,-.5) -- node[below]{$\Lambda_2^{(2)}$} (8.25,-.5);

\draw[thick,<->] (.5,.5) -- (2.75,.5);
\node[align=center, above] at (1.625,.5) {$K-R$};
\draw[thick,<->] (3,.5) -- (4.25,.5);
\node[align=center, above] at (3.625,.5) {$R$};
\draw[thick,<->] (4.5,.5) -- (6.75,.5);
\node[align=center, above] at (5.625,.5) {$K-R$};
\draw[thick,<->] (7,.5) -- (8.25,.5);
\node[align=center, above] at (7.625,.5) {$R$};

\end{tikzpicture}
\caption{Arrangement in the cell~$[1, 2K]$ for~$K=16$ and~$R=6$ }\label{f_two_scale_decomposition_of_lattice}
\end{figure}

Let us begin with decomposing the lattice with two types of blocks (cf. Figure~\ref{f_two_scale_decomposition_of_lattice}):
\begin{align}
\Lambda &: = [N] = \{1, 2, \cdots, N \}, \\
\Lambda_{1} &: = \bigcup_{ l \in \mathbb{Z}} \Lambda \cap \left( \left[1, K-R \right]  + (l-1)K \right) = \bigcup_{l =1}^{M} \Lambda_1^{(l)} \label{d_lambda1}, \\
 \Lambda_{2} &: = \bigcup_{ l \in \mathbb{Z}} \Lambda \cap  \left( \left[ K-R+1 , K \right] + (l-1)K \right) = \bigcup_{l =1}^{M} \Lambda_2 ^{(l)}, \label{d_lambda2}
\end{align}
where~$R$ is the interaction range of the particles (cf.~\eqref{e_d_hamiltonian}). The main idea is to decompose the measure~$\mu_{N,m} ( dx | y) $ into two parts as follows:
\begin{align}
\mu_{N,m} (dx|y) = \mu_{N,m} (dx^{\Lambda_1} | x^{\Lambda_2}, y) \bar{\mu}_{N,m} ( dx^{\Lambda_2} | y).
\end{align}
Again, this should be understood in a weak sense, i.e., for any test function~$\xi$,
\begin{align}
\int \xi(x) \mu_{N,m} (dx|y) = \int \left( \int \xi(x^{\Lambda_1}, x^{\Lambda_2})  \mu_{N,m} (dx^{\Lambda_1} | x^{\Lambda_2}, y) \right) \bar{\mu}_{N,m} (dx^{\Lambda_2}).
\end{align}

We prove the uniform LSI for the conditional measure~$ \mu_{N,m} (dx^{\Lambda_1} | x^{\Lambda_2}, y) $ and the marginal measure~$\bar{\mu}_{N,m} (dx^{\Lambda_2} | y) $ separately. Then uniform LSI for the full measure~$\mu_{N, m}$ is deduced via the two-scale criterion for the LSI (cf.~\cite[Theorem 3]{GrOtViWe09} or~\cite[Proposition 6]{KwMe19b}). More precisely, we have

\begin{lemma}\label{l_conditional_lsi}
The conditional measure~$ \mu_{N,m} (dx^{\Lambda_1} | x^{\Lambda_2}, y) $ satisfies LSI($\varrho_1$), where $\varrho_1>0$ is a constant independent of the system size~$|\Lambda_1|$, the external field~$s$, the mean spin~$m$, conditioned spins~$x^{\Lambda_2}$, and the macroscopic state~$y$.
\end{lemma}

\begin{lemma}\label{l_marginal_lsi} The marginal measure~$\bar{\mu}_{N,m} (dx^{\Lambda_2} | y) $ satisfies LSI($\rho_2$), where~$\rho_2>0$ is independent of the external field~$s$, the mean spin~$m$, and the macroscopic state~$y$.
\end{lemma}

\begin{lemma}\label{l_combine_lsi} Assume that
\begin{itemize}
    \item The conditional measure~$ \mu_{N,m} (dx^{\Lambda_1} | x^{\Lambda_2}, y) $ satisfies LSI($\varrho_1$), where $\varrho_1$ is a positive constant independent of~$|\Lambda_1|$,~$s$,~$m$,~$x^{\Lambda_2}$ and~$y$.
    \item The marginal measure~$\bar{\mu}_{N,m} (dx^{\Lambda_2} | y) $ satisfies LSI($\rho_2$), where~$\rho_2$ is a positive constant independent of~$s$,~$m$ and~$y$.
\end{itemize}
Then the ce~$\mu_{N,m}$ satisfies LSI($\rho$), where~$\rho$ is a positive constant independent of the system size~$N$, the external field~$s$, and the mean spin~$m$.
\end{lemma}

\medskip

\begin{itemize}
    \item Lemma~\ref{l_conditional_lsi} is a consequence of Lemma~\ref{p_uniform_lsi} and tensorization principle(cf. Theorem~\ref{a_tensorization}). Indeed, by conditioning on the spins~$x^{\Lambda_2}$, the blocks~$\Lambda_1^{(l)}$ do not interact within the Hamiltonian. Therefore the conditional measure~$\mu_{N,m}(dx^{\Lambda_1} | x^{\Lambda_2}, y)$ tensorizes on~$\bigotimes_{l=1}^{M} X_{K-R, \tilde{y}_l}$, where
\begin{align} \label{e_modified_mean_spin}
\tilde{y}_l : = \frac{Ky_l - \sum_{j \in \Lambda_2^{(l)}}x_j}{K-R}.
\end{align}
Because each tensorized measure on~$X_{K-R, \tilde{y}_l}$ has the same structure as one dimensional ce~$\mu_{N,m}$, it satisfies a uniform LSI by Lemma~\ref{p_uniform_lsi}. Then an application of Theorem~\ref{a_tensorization} yields Lemma~\ref{l_conditional_lsi}. \\

\item The proof of Lemma~\ref{l_marginal_lsi} utilizes the Otto-Reznikoff Criterion (cf. Theorem~\ref{a_otto_reznikoff}). The details are provided in Section~\ref{s_proof_marginal_lsi}. \\

\item The proof of Lemma~\ref{l_combine_lsi} is almost identical to that of~\cite[Proposition 6]{KwMe19b}. We refer to~\cite{KwMe19b} or~\cite{GrOtViWe09} for more details. \\
\end{itemize}

With the help of the lemmas above (Lemma~\ref{l_conditional_lsi}, Lemma~\ref{l_marginal_lsi} and Lemma~\ref{l_combine_lsi}) we establish the desired statement: \\

\noindent \emph{Proof of Theorem~\ref{p_generalized_lsi}.} \ A combination of Lemma~\ref{l_conditional_lsi}, Lemma~\ref{l_marginal_lsi} and Lemma~\ref{l_combine_lsi} proves Theorem~\ref{p_generalized_lsi}.
\qed

\medskip

\subsection{Proof of Lemma~\ref{l_marginal_lsi}} \label{s_proof_marginal_lsi}
The main idea of the proof of Lemma~\ref{l_marginal_lsi} is to apply the Otto-Reznikoff Criterion (Theorem~\ref{a_otto_reznikoff}). \\

For each~$l \in \{1, \cdots, M\}$ we denote (with a slight abuse of notation)
\begin{align} \label{e_block_hamiltonian}
H_l(x^{B(l)} | \bar{x}^{B(l)}) = \sum_{i \in B(l)} \psi (x_i) + \frac{1}{2}\sum_{i, j \in B(l)} M_{ij} x_i x_j + \sum_{\substack{i \in B(l) \\ k \notin B(l)}}M_{ij}x_i x_k
\end{align}
and
\begin{align}
H_{\bar{l}}(\bar{x}^{B(l)}) &= H(x) - H_l(x^{B(l)} | \bar{x}^{B(l)}) \\
& = \sum_{i \notin B(l)} \psi (x_i) + \frac{1}{2} \sum_{i, j \notin B(l)} M_{ij} x_i x_j.
\end{align}

The Hamiltonian Q associated to the marginal measure~$\bar{\mu}_{N,m} (dx^{\Lambda_2} |y )$ is
\begin{align}
Q(dx^{\Lambda_2} | y ) = -\log \int_{ \substack{ \frac{1}{K-R} \sum_{i \in \Lambda_1 ^{(l)}} x_i = \tilde{y}_l \\ l=1, \cdots, M  } }\exp\left(-H(x)\right) \mathcal{L}(dx^{\Lambda_1}),
\end{align}
where~$\tilde{y}_l$ is given by~\eqref{e_modified_mean_spin}. In particular, a rearrangement of the integral gives
\begin{align}
Q(x^{\Lambda_2} | y) = - \log \int_{\substack{ \frac{1}{K-R} \sum_{i \in \Lambda_1 ^{(k)}} x_i = \tilde{y}_k \\ k \in [M] \setminus \{l \}  }}\exp\left(-H(\bar{x}^{B(l)})\right)\exp \left( -  Q_l (x^{\Lambda_2 ^{(l)}} | \bar{x}^{B(l)} ) \right)  \mathcal{L}(d \bar{x}^{\Lambda_1^{(l)}}),
\end{align}
where~$Q_l$ is the block Hamiltonian defined by
\begin{align}
Q_l (x^{\Lambda_2 ^{(l)}} | \bar{x}^{B(l)} ) = - \log  \int_{\frac{1}{K-R}\sum_{i \in \Lambda_1^{(l)}} x_i = \tilde{y}_l  } \exp\left( - H(x^{B(l)} | \bar{x}^{B(l)} ) \right)\mathcal{L}(d x^{\Lambda_1^{(l)}}).
\end{align}

In~\cite{KwMe19b}, the authors deduced that~$Q_l$ can be decomposed into a sum of strictly convex function~$\tilde{\Psi}_l^c$ and bounded perturbation~$\tilde{\Psi}_l^b$.

\begin{lemma}[(15) and (16) in~\cite{KwMe19b}] \label{l_block_hamiltonian_decomposition}
There exist functions~$\tilde{\Psi}_l^c$ and~$\tilde{\Psi}_l^b$ such that
\begin{itemize}
    \item $Q_l = \tilde{\Psi}_l^c + \tilde{\Psi}_l^b$.
    \item For block size~$K$ large enough,~$\tilde{\Psi}_l^c$ is strictly convex.
    \item The function~$\tilde{\Psi}_l^b$ is uniformly bounded.
\end{itemize}
Moreover, the strict convexity of~$\tilde{\Psi}_l^c$ and boundedness of~$\tilde{\Psi}_l^b$ is independent of the conditioned spins~$\bar{x}^{B(l)}$.
\end{lemma}

The first step towards the proof of Lemma~\ref{l_marginal_lsi} is to prove that each block marginal measure~$\bar{\mu}_{N,m}(dx^{\Lambda_2 ^{(l)}} | y, \bar{x}^{\Lambda_2 ^{(l)}} )$ satisfies a uniform LSI.

\begin{lemma} \label{l_block_marginal_lsi}
The block marginal measure~$\bar{\mu}_{N,m}(dx^{ \Lambda_2 ^{(l)}} | y, \bar{x}^{\Lambda_2^{(l)}} )$ satisfies a uniform LSI.
\end{lemma}

\noindent \emph{Proof of Lemma~\ref{l_block_marginal_lsi}.} \ Let us fix~$l \in \{1, \cdots, M\}$ and decompose~$Q_l$ into strictly convex part~$\tilde{\Psi}_l ^c$ and bounded perturbation part~$\tilde{\Psi}_l^b$ (Lemma~\ref{l_block_hamiltonian_decomposition}). Our aim is to decompose~$Q$ into~$\tilde{\Phi}_l ^c$ and~$\tilde{\Phi}_l^b$ so that when restricted to the spins~$x^{\Lambda_2 ^{(1)}}$,~$\tilde{\Phi}_l ^c$ is strictly convex and~$\tilde{\Phi}_l^b$ is bounded. Then the desired statement follows from Bakry-\'{E}mery criterion (Theorem~\ref{a_bakry_emery}) and Holley-Stroock Perturbation Principle (Theorem~\ref{a_holley_stroock}). \\

To see this, let us decompose~$Q$ as follows:

\begin{align}
Q(dx^{\Lambda_2} | y) &= - \log \int_{\substack{ \frac{1}{K-R} \sum_{i \in \Lambda_1 ^{(k)}} x_i = \tilde{y}_k \\ k \in [M] \setminus \{l \} }}\exp\left(-H(\bar{x}^{B(l)}\right) \exp\left( -\tilde{\Psi}_l ^c - \tilde{\Psi}_l ^b \right) \mathcal{L}(d \bar{x}^{\Lambda_1^{(l)}}) \\
& = - \log \int_{\substack{ \frac{1}{K-R} \sum_{i \in \Lambda_1 ^{(k)}} x_i = \tilde{y}_k \\ k \in [M] \setminus \{l \}  }}\exp\left(-H(\bar{x}^{B(l)}\right) \exp\left( -\tilde{\Psi}_l ^c  \right) \mathcal{L}(d \bar{x}^{\Lambda_1^{(l)}}) \\
& \quad + \left( \log \int_{\substack{ \frac{1}{K-R} \sum_{i \in \Lambda_1 ^{(k)}} x_i = \tilde{y}_k \\ k \in [M] \setminus \{l \}  }}\exp\left(-H(\bar{x}^{B(l)}\right) \exp\left( -\tilde{\Psi}_l ^c  \right) \mathcal{L}(d \bar{x}^{\Lambda_1^{(l)}})  \right. \\
& \qquad \left.  - \log \int_{\substack{ \frac{1}{K-R} \sum_{i \in \Lambda_1 ^{(k)}} x_i = \tilde{y}_k \\ k \in [M] \setminus \{l \} }}\exp\left(-H(\bar{x}^{B(l)}\right) \exp\left( -\tilde{\Psi}_l ^c - \tilde{\Psi}_l ^b \right) \mathcal{L}(d \bar{x}^{\Lambda_1^{(l)}}) \right) \\
& = : \tilde{\Phi}_l ^c + \tilde{\Phi}_l ^b.
\end{align}

\textbf{Step 1.} Strict convexity of~$\tilde{\Phi}_l ^c$, restricted to the spins~$x^{\Lambda_2 ^{(1)}}$. \\

Let~$\tilde{\mu}_{l}$ be the probability measure with density
\begin{align}
\tilde{\mu}_{l} (d \bar{x}^{\Lambda_1^{(l)}}) &= \frac{1}{Z} \mathds{1}_{\substack{ \frac{1}{K-R} \sum_{i \in \Lambda_1 ^{(k)}} x_i = \tilde{y}_k \\ k \in [M] \setminus \{l \}  }} \exp\left(-H(\bar{x}^{B(l)})\right) \exp \left( -\tilde{\Psi}_l ^c  \right) \mathcal{L}(d \bar{x}^{\Lambda_1^{(l)}}) \\
& =\frac{1}{Z} \mathds{1}_{\substack{ \frac{1}{K-R} \sum_{i \in \Lambda_1 ^{(k)}} x_i = \tilde{y}_k \\ k \in [M] \setminus \{l \}  }} \exp\left(-H(\bar{x}^{B(l)})  -\tilde{\Psi}_l ^c  \right) \mathcal{L}(d \bar{x}^{\Lambda_1^{(l)}}).
\end{align}

We note that~$H(\bar{x}^{B(l)}) $ is independent of the spins~$x^{\Lambda_2 ^{(l)}}$ and thus
\begin{align}
\frac{d}{dx_i} \tilde{\Psi}_l ^c = \frac{d}{dx_i} \left( H(\bar{x}^{B(l)}) + \tilde{\Psi}_l ^c \right) \qquad \text{for } i \in \Lambda_2 ^{(l)}.
\end{align}
Therefore a straightforward calculation yields that for any $i, j \in \Lambda_2^{(l)}$,
\begin{align}
\frac{d^2}{dx_i dx_j} \tilde{\Phi}_l ^c = \mathbb{E}_{\tilde{\mu}_{N,m}} \left[ \frac{d^2}{dx_i dx_j} \tilde{\Psi}_l ^c \right] - \cov_{\tilde{\mu}_{N,m}} \left( \frac{d}{dx_i} \tilde{\Psi}_l ^c, \frac{d}{dx_j} \tilde{\Psi}_l ^c \right).
\end{align}
Due to the strict convexity of~$\tilde{\Psi}_l ^c$, an application of Brascamp-lieb inequality implies that~$\tilde{\Phi}_l ^c$ is also uniformly strictly convex on $\mathbb{R}^{\Lambda_2 ^{(l)}}$.  \\

\textbf{Step 2.} Boundedness of~$\tilde{\Phi}_l ^b$, restricted to the spins~$x^{\Lambda_2 ^{(1)}}$. \\

We write
\begin{align}
     \tilde{\Phi}_l ^b &=  \log \int_{\substack{ \frac{1}{K-R} \sum_{i \in \Lambda_1 ^{(k)}} x_i = \tilde{y}_k \\ k \in [M] \setminus \{l \}  }}\exp\left(-H(\bar{x}^{B(l)}\right) \exp\left( -\tilde{\Psi}_l ^c  \right) \mathcal{L}(d \bar{x}^{\Lambda_1^{(l)}})   \\
& \qquad   - \log \int_{\substack{ \frac{1}{K-R} \sum_{i \in \Lambda_1 ^{(k)}} x_i = \tilde{y}_k \\ k \in [M] \setminus \{l \} }}\exp\left(-H(\bar{x}^{B(l)}\right) \exp\left( -\tilde{\Psi}_l ^c - \tilde{\Psi}_l ^b \right) \mathcal{L}(d \bar{x}^{\Lambda_1^{(l)}}) \\
& = - \log \frac{ \int_{\substack{ \frac{1}{K-R} \sum_{i \in \Lambda_1 ^{(k)}} x_i = \tilde{y}_k \\ k \in [M] \setminus \{l \} }}\exp\left(-H(\bar{x}^{B(l)}\right) \exp\left( -\tilde{\Psi}_l ^c - \tilde{\Psi}_l ^b \right) \mathcal{L}(d \bar{x}^{\Lambda_1^{(l)}})} { \int_{\substack{ \frac{1}{K-R} \sum_{i \in \Lambda_1 ^{(k)}} x_i = \tilde{y}_k \\ k \in [M] \setminus \{l \}  }}\exp\left(-H(\bar{x}^{B(l)}\right) \exp\left( -\tilde{\Psi}_l ^c  \right) \mathcal{L}(d \bar{x}^{\Lambda_1^{(l)}})} \\
& = -\log \mathbb{E}_{\tilde{\mu}_l} \left[ \exp \left( - \tilde{\Psi}_l ^b \right)\right].
\end{align}
Therefore the boundedness of~$\tilde{\Phi}_l ^b$ follows from the boundedness of~$\tilde{\Psi}_l ^b$. 

\qed

Next, we shall prove that the interactions between blocks~$\Lambda_2 ^{(l)}$ and~$\Lambda_2 ^{(n)}$ become small enough when choosing the block size~$K$ large enough. Let us first introduce an auxiliary statement.

\begin{lemma}[Lemma 16 in~\cite{KwMe19b}] \label{l_derivative_change_of_var}
It holds that
\begin{align}
\frac{d^2}{dx_i dx_j} Q (dx^{\Lambda_2} | y)  = - \cov_{\mu_{N,m}(dx^{ \Lambda_1} | x^{\Lambda_2}, y)} \left( \frac{\partial}{\partial x_i} H(x)- \frac{\partial}{\partial x_{i-R}} H(x), \frac{\partial}{\partial x_j} H(x) - \frac{\partial}{\partial x_{j-R}} H(x)  \right).
\end{align}
\end{lemma}

The following statement is the second main ingredient for proving Lemma~\ref{l_marginal_lsi}.

\begin{lemma} \label{l_block_marginal_interaction}
For any~$n, l \in \{1, \cdots, M \}$ with~$n \ne l$, it holds that
\begin{align}
\left|\frac{d^2}{dx_i dx_j} Q (dx^{\Lambda_2} | y) \right| \lesssim \frac{R}{K} + R \exp \left( -CK |n-l| \right) \qquad \text{for all } i \in \Lambda_2 ^{(n)}, j \in \Lambda_2 ^{(l)}.
\end{align}
\end{lemma}

\medskip

\noindent \emph{Proof of Lemma~\ref{l_block_marginal_interaction}.} \ Let us begin with a simple observation. Let~$f$ and~$g$ be functions supported on~$\Lambda_1 ^{(l)}$ and~$\Lambda_1 ^{(n)}$, respectively. Because the measure~$\mu_{N,m}(dx^{\Lambda_1}|x^{\Lambda_2}, y)$ tensorizes on~$\bigotimes_{l=1}^{M} X_{K-R, \tilde{y}_l}$,  it holds that
\begin{align} \label{e_tensor_covariance}
\cov_{\mu_{N,m}(dx^{\Lambda_1}| x^{\Lambda_2}, y)} (f, g) =
\begin{cases} 0 \qquad \qquad \qquad \qquad \qquad &\text{if } l \neq n \\
\cov_{\mu_{N,m}(dx^{\Lambda_1 ^{(l)}} | x^{\Lambda_2}, y)} (f, g) \qquad &\text{otherwise}
\end{cases}
\end{align}
Then a combination of Lemma~\ref{l_derivative_change_of_var},~\eqref{e_tensor_covariance} and Proposition~\ref{p_decay_of_correlations_ce} yields
\begin{align}
\left|\frac{d^2}{dx_i dx_j} \bar{Q} (dx^{\Lambda_2} | y) \right| \lesssim \frac{R}{K} + R \exp \left( -CK |n-l| \right).
\end{align}
\qed

Now we are ready to present the proof of Lemma~\ref{l_marginal_lsi}. \\

\noindent \emph{Proof of Lemma~\ref{l_marginal_lsi}.} \ By choosing the block size~$K$ large enough, the Otto-Reznikoff Criterion (Theorem~\ref{a_otto_reznikoff}) applied with the help of Lemma~\ref{l_block_marginal_lsi} and Lemma~\ref{l_block_marginal_interaction} finishes the proof of Lemma~\ref{l_marginal_lsi}.
\qed

\section{Strict convexity of coarse-grained Hamiltonian: proof of Theorem~\ref{p_strict_convexity_cg_hamiltonian}} \label{s_proof_strict_convexity_cg_hamiltonian}

The proof of Theorem~\ref{p_strict_convexity_cg_hamiltonian} consists of two ingredients. The first one is uniform estimates of the diagonals of~$\Hess_Y \bar{H} (y)$ and the second one is the control of off-diagonal terms of~$\Hess_Y \bar{H} (y)$. For the diagonals, we adapt the one dimensional result. That is, one-dimensional coarse-grained Hamiltonian is uniformly strictly convex (cf.~\cite{KwMe18}). The off-diagonal terms are controlled by decay of correlations. (cf. Theorem~\ref{p_decay_of_correlations_ce}) \\

\begin{lemma}\label{l_strict_convexity_diagonal} There is a positive constant~$\tau>0$ such that for each~$l \in \{1, \cdots, M\}$, it holds that
\begin{align}
0< \tau \leq \left( \Hess_Y \bar{H} (y)\right)_{ll} \leq \frac{1}{\tau} < \infty.
\end{align}
\end{lemma}
\medskip

\begin{lemma}\label{l_smallness_off_diagonal} For each~$1 \leq l \neq n \leq M$, it holds that
\begin{align}
\left| \left( \Hess_Y \bar{H} (y) \right)_{ln}\right| \lesssim \frac{1}{K}.
\end{align}
\end{lemma}
\medskip

The proof of Lemma~\ref{l_strict_convexity_diagonal} and Lemma~\ref{l_smallness_off_diagonal} are given in Section~\ref{s_proof_strict_convexity_diagonal} and Section~\ref{s_proof_smallness_off_diagonal}, respectively. Let us now proceed to the proof of Proposition~\ref{p_strict_convexity_cg_hamiltonian}. \\

\noindent \emph{Proof of Theorem~\ref{p_strict_convexity_cg_hamiltonian}.} \ It follows directly from Lemma~\ref{l_strict_convexity_diagonal} and Lemma~\ref{l_smallness_off_diagonal} by choosing~$K$ large enough. \qed

\bigskip 

\subsection{Proof of Lemma~\ref{l_strict_convexity_diagonal}} \label{s_proof_strict_convexity_diagonal} 

Recall the definition~\eqref{e_block_hamiltonian} of~$H( x^{B(l)} | \bar{x}^ {B(l)} )$. The coarse-grained Hamiltonian associated to~$H(x^{B(l)} | \bar{x}^{B(l)})$ is
\begin{align}
\bar{H}(y_l | \bar{x}^{B(l)} )= -\frac{1}{K} \log \int_{\frac{1}{K}\sum_{i \in B(l)}x_i = y_l} \exp \left(- H(x^{B(l)} | \bar{x}^{B(l)})\right) \mathcal{L}^{K-1}(dx^{B(l)}).
\end{align}
Observing that the Hamiltonian~$H(x^{B(l)} | \bar{x}^{B(l)})$ has the same structure as the Hamiltonian given by~\eqref{e_d_hamiltonian}, a straightforward calculation using Lemma~\ref{l_hessian_cg_hamiltonian} yields the following statement:

\begin{lemma}[(34) in~\cite{Me11}] \label{l_derivative_formula}
For any~$l \in \{ 1, \cdots M \}$, it holds that
\begin{align}
&\left(\Hess_Y \bar{H} (y) \right)_{ll} \\
&= \int \frac{d^2}{dy_l ^2} \bar{H}( y_l | \bar{x}^{B(l)} ) \bar{\mu}_{N,m}(d \bar{x}^{B(l)} | y)  \\
& \quad  - \frac{1}{K} \var_{\bar{\mu}_{N,m}(d\bar{x}^{B(l)} | y)} \left( \int \left( \sum_{j \in B(l)} \left(  \sum_{i=1}^N M_{ij}x_i \right) + \psi_b ' (x_j) \right) \mu_{N,m}(dx^{B(l)} | \bar{x}^{B(l)} , y) \right). \label{e_derivative_formula}
\end{align}
\end{lemma}
\medskip

In~\cite{KwMe18}, the authors proved that under the assumption of lower bound of variance of the mean spin of the modified gce~$\mu_N^{\sigma}$, the one-dimensional coarse-grained Hamiltonian is uniformly strictly convex (cf.~\cite[Corollary 2]{KwMe18}). Combined with Lemma~\ref{l_variance_estimate} we get 

\begin{lemma} [Extension of Corollary 2 in~\cite{KwMe18}] \label{l_hessian_diagonal}
There is a positive constant~$\lambda$ independent of the system size~$N$, the external field~$s$, and the mean spin~$m$ such that
\begin{align} \label{e_hessian_diagonal}
\lambda \leq \frac{d^2}{dy_l ^2} \bar{H} ( y_l | \bar{x}^{B(l)} ) \leq \frac{1}{\lambda}.
\end{align}
\end{lemma}

The next statement implies that the right hand side of~\eqref{e_derivative_formula} is small when~$K$ is large enough.

\begin{lemma} \label{l_hessian_off_diagonal}
It holds that
\begin{align} \label{e_off_diagonal}
\frac{1}{K} \var_{\bar{\mu}_{N,m}(d\bar{x}^{B(l)} | y)} \left( \int \left( \sum_{j \in B(l)} \left(  \sum_{i=1}^N M_{ij}X_i \right) + \psi_b ' (X_j) \right) \mu_{N,m}(dx^{B(l)} | \bar{x}^{B(l)} , y) \right) \lesssim \frac{1}{K}
\end{align}
\end{lemma}

\noindent \emph{Proof of Lemma~\ref{l_hessian_off_diagonal}.} \ Note that the conditional measure~$\mu_{N,m}(dx^{B(l)} | \bar{x}^{B(l)} , y) $ is given by
\begin{align}
\mu_{N,m}(dx^{B(l)} | \bar{x}^{B(l)} , y) & = \frac{1}{Z} \mathds{1}_{ \left\{ \frac{1}{K} \sum_{i \in B(l)} x_i =y_l \right\}}(x^{B(l)}) \exp \left(- H(x^{B(l)} | \bar{x}^{B(l)}) \right) \mathcal{L}^{K-1}(dx^{B(l)}).
\end{align}
For each~$l \in \{1, \cdots, M\}$ define~$E_l : = \{k \notin B(l) : \exists i \in B(l) \text{ such that } |i-k| \leq R  \}$. By cancellation of constant terms with the partition function, the term
\begin{align}
\int \left( \sum_{j \in B(l)} \left(  \sum_{i=1}^N M_{ij}x_i \right) + \psi_b ' (x_j) \right) \mu_{N,m}(dx^{B(l)} | \bar{x}^{B(l)} , y)  \label{e_inside_var}
\end{align}
depends only on the spins~$x_k$, $k \in E_l$. In particular,~\eqref{e_inside_var} is a function of~$\bar{x}^{B(l)}$. This implies
\begin{align}
&\var_{\bar{\mu}_{N,m}(d\bar{x}^{B(l)} | y)} \left( \int \left( \sum_{j \in B(l)} \left(  \sum_{i=1}^N M_{ij}x_i \right) + \psi_b ' (x_j) \right) \mu_{N,m}(dx^{B(l)} | \bar{x}^{B(l)} , y) \right)  \\
& \qquad  = \var_{\mu_{N,m}(dx|y)} \left( \int \left( \sum_{j \in B(l)} \left(  \sum_{i=1}^N M_{ij}x_i \right) + \psi_b ' (x_j) \right) \mu_{N,m}(dx^{B(l)} | \bar{x}^{B(l)} , y) \right).
\end{align}
Then an application of Poincar\'e inequality for~$\mu_{N,m}(dx|y)$ (cf. Proposition~\ref{p_generalized_lsi}) yields
\begin{align}
&\var_{\mu_{N,m}(dx|y)} \left( \int \left( \sum_{j \in B(l)} \left(  \sum_{i=1}^N M_{ij}x_i \right) + \psi_b ' (x_j) \right) \mu_{N,m}(dx^{B(l)} | \bar{x}^{B(l)} , y) \right) \\
&\qquad  \lesssim \int \left| \nabla \left( \int \left( \sum_{j \in B(l)} \left(  \sum_{i=1}^N M_{ij}x_i \right) + \psi_b ' (x_j) \right) \mu_{N,m}(dx^{B(l)} | \bar{x}^{B(l)} , y) \right) \right|^2 \mu_{N,m}(dx|y).
\end{align}
Because~\eqref{e_inside_var} depends only on the spins~$x_k$,~$k \in E_l$, it holds that
\begin{align}
&\left| \nabla \left( \int \left( \sum_{j \in B(l)} \left(  \sum_{i=1}^N M_{ij}x_i \right) + \psi_b ' (x_j) \right) \mu_{N,m}(dx^{B(l)} | \bar{x}^{B(l)} , y) \right) \right| \\
&\qquad = \sum_{k \in E_l} \left( \frac{\partial}{\partial x_k} \int \left( \sum_{j \in B(l)} \left(  \sum_{i=1}^N M_{ij}x_i \right) + \psi_b ' (x_j) \right) \mu_{N,m}(dx^{B(l)} | \bar{x}^{B(l)} , y) \right)^2 \\
& \qquad = \sum_{k \in E_l} \left( \sum_{j \in B(l)} M_{kj} - \cov_{\mu_{N,m}(dx^{B(l)} | \bar{x}^{B(l)}, y)} \left( \sum_{j \in B(l)} \left(  \sum_{i=1}^N M_{ij}x_i \right) + \psi_b ' (x_j),   \sum_{\substack{i \in B(l) \\ |k-i| \leq R }}M_{ik}x_i  \right)  \right)^2. \label{e_gradient_formula}
\end{align}
Lastly, an application of Proposition~\ref{p_decay_of_correlations_ce} combined with the fact that~$|E_l| \leq 2R$ implies 
\begin{align}
T_{\eqref{e_gradient_formula}} \lesssim 1.
\end{align}
This finishes the proof of Lemma~\ref{l_hessian_off_diagonal}. \qed
\medskip

\noindent \emph{Proof of Lemma~\ref{l_strict_convexity_diagonal}.} \ Lemma~\ref{l_strict_convexity_diagonal} is a direct consequence of Lemma~\ref{l_derivative_formula}, Lemma~\ref{l_hessian_diagonal}, and Lemma~\ref{l_hessian_off_diagonal}. Indeed, by choosing~$K$ large enough, we can find a positive constant~$\tau >0$ such that
\begin{align}
0 < \tau \leq \left( \Hess_Y \bar{H} (y) \right)_{ll} \leq \frac{1}{\tau}  < \infty.
\end{align}
\qed 
\bigskip 

\subsection{Proof of Lemma~\ref{l_smallness_off_diagonal}.} \label{s_proof_smallness_off_diagonal} 

The main ingredient for the proof of Lemma~\ref{l_smallness_off_diagonal} is the following representation of the Hessian of~$\bar{H}$
\begin{lemma}[Lemma 2 in~\cite{Me11}] \label{l_hessian_formula}
For any~$1 \leq l \neq  n \leq M$, it holds that
\begin{align}
&\left( \Hess_Y \bar{H}(y) \right)_{ln}\\
&= \frac{1}{K}\sum_{i \in B(l), j \in B(n)} M_{ij} \label{e_interaction_hessian}\\
& \quad - \frac{1}{K} \cov_{\mu_{N,m}(dx|y)} \left( \sum_{j \in B(l)} \left( \sum_{i=1}^N M_{ij}x_i\right) + \psi_b ' (x_j) , \sum_{j \in B(n)} \left( \sum_{i=1}^N M_{ij}x_i\right) + \psi_b ' (x_j)  \right). \label{e_covariance_hessian}
\end{align}
\end{lemma}
\medskip 

\noindent \emph{Proof of Lemma~\ref{l_smallness_off_diagonal}.} \ 
We observe that for~$l \neq n$, there are at most~$R^2$ many pairs~$(i,j)$ with $i \in B(l)$,~$j \in B(n)$ and~$|i-j| \leq R$. For such~$(i,j)$, we know~$|M_{ij}|$ is uniformly bounded by 1 and hence
\begin{align}
\left|T_{\eqref{e_interaction_hessian}}\right| \lesssim R^2 \cdot \frac{1}{K} \lesssim \frac{1}{K}.
\end{align}

Let us turn to the estimation of~\eqref{e_covariance_hessian}. The law of total variance yields
\begin{align}
&\cov_{\mu_{N,m}(dx|y)} \left( \sum_{j \in B(l)} \left( \sum_{i=1}^N M_{ij}x_i\right) + \psi_b ' (x_j) , \sum_{j \in B(n)} \left( \sum_{i=1}^N M_{ij}x_i\right) + \psi_b ' (x_j)  \right) \\
& = \cov_{\bar{\mu}_{N,m}(d\bar{x}^{B(l)}|y)} \left( \int \sum_{j \in B(l)} \left( \sum_{i=1}^N M_{ij}x_i\right) + \psi_b ' (x_j) \mu_{N,m}(dx^{B(l)}|\bar{x}^{B(l)}, y) , \right. \\
&  \qquad \qquad \qquad \qquad \qquad \left. \int \sum_{j \in B(n)} \left( \sum_{i=1}^N M_{ij}x_i\right) + \psi_b ' (x_j) \mu_{N,m}(dx^{B(l)}|\bar{x}^{B(l)}, y) \right) \\
& \quad + \int \cov_{\mu_{N,m}(dx^{B(l)}|\bar{x}^{B(l)}, y)} \left( \sum_{j \in B(l)} \left( \sum_{i=1}^N M_{ij}x_i\right) + \psi_b ' (x_j) ,  \right. \\
&\qquad \qquad \qquad \qquad \qquad \qquad \qquad \qquad  \left.  \sum_{j \in B(n)} \left( \sum_{i=1}^N M_{ij}x_i\right) + \psi_b ' (x_j) \right) \bar{\mu}_{N,m}(d\bar{x}^{B(l)}|y).
\end{align}
Then an application of Proposition~\ref{p_decay_of_correlations_ce} as in Lemma~\ref{l_strict_convexity_diagonal} yields the desired estimate
\begin{align}
\left|T_{\eqref{e_covariance_hessian}} \right| \lesssim \frac{1}{K}.
\end{align}
Hence we conclude
\begin{align}
\left| \left( \Hess_Y \bar{H} (y) \right)_{ln}\right| \leq \left|T_{\eqref{e_interaction_hessian}} \right| + \left|T_{\eqref{e_covariance_hessian}}\right| \lesssim \frac{1}{K}.
\end{align}
\qed

\section{Hydrodynamic limit: proof of Theorem~\ref{p_hydrodynamic_limit}} \label{s_proof_hydro}

In this section, we provide the proof of our main theorem: hydrodynamic limit of Kawasaki dynamics (Theorem~\ref{p_hydrodynamic_limit}). The main idea for the proof of Theorem~\ref{p_hydrodynamic_limit} is the two-scale approach (cf.~\cite[Theorem 8]{GrOtViWe09}). In~\cite{GrOtViWe09}, the hydrodynamic limit of the Kawasaki dynamics was deduced via two-scale approach where there is no-interactions within the Hamiltonian (see~\cite[Theorem 17]{GrOtViWe09}). The problem becomes a lot more subtle when we add strong finite-range interactions within the Hamiltonian. For example, because neighboring blocks interact with each other, the coarse-grained Hamiltonian~$\bar{H}_Y$ (cf.~\eqref{e_def_cg_hamiltonian}) cannot be decomposed into a sum of one-dimensional coarse-grained Hamiltonian of the form~\eqref{e_def_1d_cg_ham}. \\

We overcome this difficulty by introducing an auxiliary Hamiltonian~$H_{\text{aux}}$ obtained by removing the interactions between different blocks in the Hamiltonian~$H$ (see~\eqref{e_def_aux_hamiltonian}). Due to the finite range interactions, the amount of interactions that we remove from the formal Hamiltonian~$H$ is very small compared to the whole system size. Therefore one can expect that~$H$ and~$H_{\text{aux}}$ are \textit{close}. This allows us to take advantage of nice structure of~$H_{\text{aux}}$ (e.g. block decomposition). \\

Let us begin with introducing auxiliary definitions. First of all, we define the coarse-grained operator~$\bar{A} : Y \to Y $ by
\begin{align}
(\bar{A})^{-1} = PA^{-1} N P^*.
\end{align}
For given~$\eta_0 \in Y$, consider the mesoscopic  analog of Kawasaki dynamics:
\begin{align} \label{e_mesoscopic_parabolic}
\begin{cases}
\frac{d \eta}{dt} = - \bar{A} \nabla_Y \bar{H}_Y(\eta ) \\
\eta (0) = \eta_0.
\end{cases}
\end{align}

Recalling the identification of~$X$ and~$\bar{X}$ (see Section~\ref{s_main_hydro}), we identify~$Y$ with the space~$\bar{Y}$ of piecewise constant, mean~$m$ functions on~$\mathbb{T}^1 = \mathbb{R} \backslash \mathbb{Z}$ defined by
\begin{align} \label{e_def_barX}
\bar{Y} : = \left\{ \bar{y} : \mathbb{T}^1 \to \mathbb{R} ; \ \bar{y} \text{ is constant on } \left( \frac{l-1}{M}, \frac{l}{M} \right] \text{ for } l=1, \cdots, M, \text{ and has mean } m \right\}.
\end{align}

\medskip

The main idea of the two-scale approach is to prove the closeness of microscopic-mesoscopic solutions and mesoscopic-macroscopic solutions. \\

Consider a sequence~$\{M_\nu, N_\nu\}_{\nu=1}^{\infty}$ such that
\begin{align}
M_\nu \uparrow \infty, \qquad N_\nu \uparrow \infty, \qquad K_\nu =\frac{N_\nu}{M_\nu} \uparrow \infty.
\end{align}
This means that the size of each block and the number of blocks are simultaneously increasing to the  infinity.

\begin{convention}
Following the convention of~\cite{GrOtViWe09}, we write~$M, N, K$ for~$M_{\nu}, N_{\nu}, K_{\nu}$. We also denote~$X = X_{N^{\nu}, m}$,~$Y = Y_{M^{\nu}, m}$, and so on in the remaining sections.
\end{convention}

For given~$\zeta_0$, choose a sequence of step functions~$\{ \bar{\eta}_0 ^{\nu} \}_{\nu =1}^{\infty}$ in~$\bar{Y}$ that converges to~$\zeta_0$ in~$L^2$:
\begin{align}
    \| \bar{\eta}_0 ^{\nu} - \zeta_0 \|_{L^2} \to 0 \qquad \text{as } \nu \to \infty. \label{d_eta0}
\end{align}

For each~$\nu$, let~$\eta_0 ^{\nu} \in Y$ be the vector that corresponds to the step function~$\bar{\eta}_0 ^{\nu}$, and denote ~$\eta^{\nu}$ by a  solution of the mesoscopic parabolic equation \eqref{e_mesoscopic_parabolic} with the initial data $\eta (0) = \eta_0 ^{\nu}$.


The first main ingredient is the closeness of microscopic-mesoscopic solutions.

\begin{proposition} \label{p_micro_meso}
For any~$T >0$, it holds that
\begin{align}
\lim_{\nu \to \infty} \sup_{0 \leq t \leq T} \int  \| \bar{x} - \bar{\eta} ^{\nu}(t, \cdot) \|_{H^{-1}} f(t, x) \mu_{N, m}(dx) = 0.
\end{align}
\end{proposition}

The second ingredient is the closeness of mesoscopic-macroscopic solutions. Note that by the expression \eqref{211} and Proposition \ref{l_convergence_cg_ham}, a function $\varphi$ in the macroscopic parabolic equation \eqref{e_heat_eqn} is a pointwise limit of the one-dimensional coarse-grained Hamiltonian $\bar{H}_N$ and is differentiable.

\begin{proposition} \label{p_meso_macro}
The step functions~$\bar{\eta}^{\nu}$ converge strongly in~$L^{\infty}(H^{-1})$ to the unique weak solution~$\zeta$ of~\eqref{e_heat_eqn}. In particular, it holds that for any~$T >0$,
\begin{align}
\lim_{\nu \to \infty} \sup_{0 \leq t \leq T} \| \bar{\eta}^{\nu} (t, \cdot) - \zeta(t, \cdot) \|_{H^{-1}}^2 = 0
\end{align}
\end{proposition}

We provide the proof of Proposition~\ref{p_micro_meso} and Proposition~\ref{p_meso_macro} in Section~\ref{s_proof_micro_meso} and Section~\ref{s_proof_meso_macro}, respectively. \\

\noindent \emph{Proof of Theorem~\ref{p_hydrodynamic_limit}.} \ Following the notations from above, Proposition~\ref{p_micro_meso} and Proposition~\ref{p_meso_macro} imply
\begin{align}
&\lim_{N \to \infty} \sup_{0 \leq t \leq T} \int \| \bar{x} - \zeta (t, \cdot) \|_{H^{-1}}^2 f(t,x) \mu(dx)\\
&\qquad \leq 2 \lim_{N \to \infty} \sup_{0 \leq t \leq T} \left(  \int \| \bar{x} - \bar{\eta}^{v} (t, \cdot) \|_{H^{-1}}^2 f(t,x) \mu(dx)  + \int \| \bar{\eta}^{v} - \zeta (t, \cdot) \|_{H^{-1}}^2 f(t,x) \mu(dx) \right) \\
&\qquad  = 2 \lim_{N \to \infty} \sup_{0 \leq t \leq T} \left(  \int \| \bar{x} - \bar{\eta}^{v} (t, \cdot) \|_{H^{-1}}^2 f(t,x) \mu(dx)  +  \| \bar{\eta}^{v} - \zeta (t, \cdot) \|_{H^{-1}}^2  \right) \\
&\qquad  = 0. 
\end{align}
\qed

\subsection{Proof of Proposition~\ref{p_micro_meso}} \label{s_proof_micro_meso}
The key ingredient of the proof  of Proposition~\ref{p_micro_meso} is a two-scale criterion for the hydrodynamic limit which was originally obtained in  \cite[Theorem 8]{GrOtViWe09}. This was successfully used to establish a hydrodynamic limit of Kawasaki dynamics without interactions. We first introduce a general two-scale criterion for the hydrodynamic limit developed in \cite{GrOtViWe09}, and then apply this to the general ce using the results established so far.

\begin{theorem}[Two-Scale Criterion for the hydrodynamic limit~\cite{GrOtViWe09}] \label{a_two_scale_hydrodynamic}
Let
\begin{align}
    \mu(dx)= \frac{1}{Z}\exp\left(-H(x)\right)dx
\end{align}
be a probability measure on~$X$. Assume a linear operator~$P : X \to Y$ satisfies~$PNP^*= \textup{Id}_Y$ for some large~$N \in \mathbb{N}$. Assume further the following:

(\romannumeral 1). It holds that 
\begin{align}
\kappa : = \max \left \{ \langle \Hess H(x) \cdot u, v \rangle : \ u \in \textup{Ran}(NP^* P), v \in \textup{Ran}(\Id_X - NP^* P), |u|=|v|=1    \right \} < \infty.
\end{align}

(\romannumeral 2). There is~$\rho>0$ such that~$\mu(dx|y)$ satisfies LSI$(\rho)$ for all~$y$

(\romannumeral 3). There is~$\lambda>0$ such that~$\langle \tilde{y}, \Hess \bar{H}(y) \tilde{y} \rangle_Y \geq \lambda \langle \tilde{y}, \tilde{y} \rangle_Y$

(\romannumeral 4). There is~$\alpha >0 $ such that~$ \int |x|^2 \mu(dx) \leq \alpha N$

(\romannumeral 5). There is~$\beta >0$ such that~$\inf_{y \in Y} \bar{H}(y) \geq - \beta$

Define~$M:= \textup{dim} Y$ and let~$A: X \to X$ be a symmetric linear operator such that:

(\romannumeral 6). There is~$\gamma>0$ such that for all~$x \in X$, $|(\textup{Id}_X -NP^* P)x|^2 \leq \gamma M^{-2}\langle x, Ax \rangle_X$

Let~$f(t,x)$ and~$\eta (t)$ solve~\eqref{e_kawasaki_de} and
\begin{align}
\frac{d \eta}{dt} = - \bar{A} \nabla_Y \bar{H}(\eta),
\end{align}
with initial data~$f(0,\cdot)$, and~$\eta_0$, where~$(\bar{A})^{-1}= PA^{-1}NP^*$. Assume

(\romannumeral 7). $\int f(0,x) \log f(0,x) \mu(dx) \leq C_1 N$, $\bar{H}(\eta_0) \leq C_2$.

Define
\begin{align}
\Theta(t): = \frac{1}{2N} \int \left\langle \left(x-NP^* \eta(t)\right),  A^{-1}\left( x- NP^* \eta (t) \right)   \right\rangle f(t,x) \mu(dx).
\end{align}
Then for any~$T>0$ we have
\begin{align}
&\max \left\{ \sup_{0 < t \leq T} \Theta (t), \frac{\lambda}{2} \int_{0}^T \left( \int_Y |y- \eta(t)|_Y ^2 \bar{f}(t,y) \bar{\mu}(dy) \right)dt  \right\} \\
& \qquad \leq \Theta (0) + T \left(\frac{M}{N}\right) + \left( \frac{C_1 \gamma \kappa^2}{2 \lambda \rho^2} \frac{1}{M^2} \right) + \left[ \sqrt{2T\gamma} \left( \alpha + \frac{2C_1}{\hat{\rho}} \right)^{\frac{1}{2}} (C_1 ^{\frac{1}{2}} + (C_2 + \beta)^{\frac{1}{2}}) \right] \frac{1}{M}.
\end{align}

\end{theorem}

As a corollary, we have both microscopic and mesoscopic closeness between the microscopic Kawasaki dynamics and the evolution \eqref{e_mesoscopic_parabolic}.
\begin{corollary}[Propagation of hydrodynamic behavior~\cite{GrOtViWe09}] \label{a_propagation}
Consider a sequence $\{X_{v}, Y_{v}, P_{v}, A_{v}, \mu_{v}, f_{0, v}, \eta_{0, v}    \}_{v=1}^{\infty}  $ of data satisfying the assumptions of Theorem~\ref{a_two_scale_hydrodynamic} for every~$\nu$ with uniform constants $\lambda, \rho, \kappa, \alpha, \beta, \gamma, C_1, C_2$. Suppose that
\begin{align}
M_v \uparrow \infty, \qquad N_v \uparrow \infty, \qquad \frac{N_v}{M_v} \uparrow \infty.
\end{align}
Further assume that
\begin{align} \label{e_cor_assumption}
\lim_{v \uparrow \infty} \frac{1}{N_v} \int (x- N_v P_v^t \eta_{0, v}) \cdot A_v ^{-1} ( x-N_v P_v^t \eta_{0, v}) f_{0,v}(x) \mu_v (dx) =0.
\end{align}
Then for any~$T>0$,
\begin{align}
\lim_{v \uparrow \infty} \sup_{0\leq t \leq T} \frac{1}{N_v} \int ( x- N_v P_v ^{t} \eta) \cdot A_v^{-1}( x- N_v P_v ^{t} \eta) f(t,x)\mu(dx) = 0,
\end{align}
and
\begin{align}
\lim_{v \uparrow \infty} \int_{0}^{T} \int_Y |y- \eta(t)|_Y ^2 \bar{f}(y) \bar{\mu}(dy) dt =0.
\end{align}
\end{corollary}

\subsubsection{Auxiliary Hamiltonian}
As mentioned at the beginning of Section  \ref{s_proof_hydro}, we introduce 
 the auxiliary Hamiltonian and related notions. First of all, define the auxiliary Hamiltonian~$H_{N, \text{aux}} = H_{\text{aux}}$ as
\begin{align}
H_{\text{aux}}(x) & := H(x) - \frac{1}{2} \sum_{l=1}^{M} \sum_{n \neq l} \sum_{\substack{i \in B(l) \\ j \in B(n)}} M_{ij} x_i x_j.  \label{e_def_aux_hamiltonian}
\end{align}
Because the interactions between different blocks are removed,~$H_{\text{aux}}$ is decomposed as follows:
\begin{align}
 H_{\text{aux}}(x)   & = \sum_{l=1}^{M}  \left( \sum_{i \in B(l)}  \left( \psi (x_i) + s_i x_i + \frac{1}{2} \sum_{\substack{j \in B(l), \\ 1 \leq |j-i| \leq R  }} M_{ij} x_i x_j \right)   \right) \\
 & = \sum_{l=1}^{M} H_K (x^{B(l)}).
\end{align}
Here, we note that there are at most~$2R^2 M$ many pairs of~$(i, j, l, n)$ such that
\begin{itemize}
    \item $l, n \in [M]$ and~$l \ne n$ 
    \item $i \in B(l)$,~$j \in B(n)$ and~$|i-j| \leq R$.
\end{itemize}

Next, the corresponding canonical ensemble~$\mu_{N, m, \text{aux}}$ is
\begin{align} 
\mu_{N, m, \text{aux}}  (dx) : &= \frac{1}{Z} \mathds{1}_{ \left\{ \frac{1}{N} \sum_{i =1}^{N} x_i =m \right\}}\left(x\right) \exp\left( - H_{\text{aux}}(x) \right) \mathcal{L}^{N-1}(dx), \label{d_ce_aux}
\end{align}
where~$\mathcal{L}^{N-1}(dx)$ denotes the~$(N-1)$-dimensional Hausdorff measure supported on~$X$.\\

Then we decompose the ce~$\mu_{N, m, \text{aux}}$ into the conditional measure~$\mu_{N, m, \text{aux}} (dx |y)$ and the marginal measure~$\bar{\mu}_{N, m, \text{aux}}(y)$ and define the corresponding coarse-grained Hamiltonian~$\bar{H}_{Y, \text{aux}}$ by
\begin{align}
    \bar{H}_{Y, \text{aux}}(y) :&= -\frac{1}{N} \log \int_{Px=y} \exp \left(-H_{\text{aux}}(x)\right) \mathcal{L}^{N-M}(dx) \\
    & = \frac{1}{M} \sum_{l=1}^{M} \left( - \frac{1}{K} \log \int_{ \{\frac{1}{N}\sum_{i \in B(l)} x_i = y_l \}}  \exp \left( -H_K (x^{B(l)}) \mathcal{L}^{M-1}(dx^{B(l)})  \right)   \right) \\
    & = \frac{1}{M} \sum_{l=1}^{M} \bar{H}_{K}  (y_l). \label{d_aux_cg_ham}
\end{align}

\begin{convention}
In Section~\ref{s_proof_micro_meso} and~\ref{s_proof_meso_macro}, we write~$\mu = \mu_{N,m}$,~$\mu_{\text{aux}} = \mu_{N, m, \text{aux}}$,~$\bar{H} = \bar{H}_Y$, and~$\bar{H}_{\text{aux}} = \bar{H}_{Y, \text{aux}}$ to reduce our notational burden.
\end{convention}

\subsubsection{Auxiliary Lemmas} \label{s_micro_meso_auxiliary}
In this section, we provide auxiliary statements that are needed in the proof of Proposition~\ref{p_micro_meso}. 
Let us begin with investigating the relationship between~$\bar{H}$ and~$\bar{H}_{\text{aux}}$. We prove that the difference between~$\bar{H}$ and~$\bar{H}_{\text{aux}}$ is small.

\begin{lemma} \label{l_close_cg_hamiltonians}
There exists a constant $C>0$ such that for any $y\in Y$,
\begin{align}
\left| \bar{H}(y)  -\bar{H}_{\text{aux}} (y) \right| \leq \frac{C}{K} \left( 1 + \|y \|_{L^2(Y)} ^2 \right).
\end{align}

\end{lemma}

\noindent \emph{Proof of Lemma~\ref{l_close_cg_hamiltonians}.} \ As there is a small cosmetic difference between the proof of Lemma~\ref{l_close_cg_hamiltonians} and Lemma~\ref{l_subadditive_gce}, we shall only outline the proof of Lemma~\ref{l_close_cg_hamiltonians}. \\

Recalling the definition~\eqref{e_def_aux_hamiltonian} of~$H_{\text{aux}}$, we have
\begin{align}
 \bar{H}_{\text{aux}}(y) - \bar{H}(y)  &=  \frac{1}{N} \log \frac{\int_{Px=y} \exp\left(-H(x)\right) \mathcal{L}^{N-M}(dx)}{\int_{Px=y} \exp\left(-H_{\text{aux}}(x)\right) \mathcal{L}^{N-M}(dx)}  \\
& = \frac{1}{N} \log \left( \mathbb{E}_{\mu_{\text{aux}}(dx|y)} \left[ - \frac{1}{2} \sum_{l=1}^{M} \sum_{n \neq l} \sum_{\substack{i \in B(l) \\ j \in B(n)}} M_{ij} x_i x_j \right] \right) \\
& \leq \frac{1}{N} \log \left( \mathbb{E}_{\mu_{\text{aux}}(dx|y)} \left[ \frac{1}{4} \sum_{l=1}^{M} \sum_{n \neq l} \sum_{\substack{i \in B(l) \\ j \in B(n)}} |M_{ij}| (x_i ^2 + x_j ^2) \right]  \right) \\
& = \frac{1}{N} \log \left( \mathbb{E}_{\mu_{\text{aux}}(dx|y)} \left[ \frac{1}{2} \sum_{l=1}^{M} \sum_{n \neq l} \sum_{\substack{i \in B(l) \\ j \in B(n)}} |M_{ij}| x_i ^2 \right]  \right) \label{e_difference_cg_hamiltonian}
\end{align}
Because the conditional measure~$\mu_{\text{aux}}(dx|y)$ tensorizes, i.e.,
\begin{align}
\mu_{\text{aux}}(dx|y) &= \bigotimes_{l=1}^{M} \frac{1}{Z} \mathds{1}_{ \left\{ \frac{1}{K} \sum_{i \in B(l)} x_i =y_l \right\}}\left(x^{B(l)}\right) \exp\left( - H_K(x^{B(l)}) \right) \mathcal{L}^{K-1}(dx^{B(l)}) \\
& =: \bigotimes_{l=1}^{M} \mu_K (dx^{B(l)} | y_l),
\end{align}
we have
\begin{align}
    T_{\eqref{e_difference_cg_hamiltonian}} =  \frac{1}{N} \sum_{l=1}^{M}  \log \left( \mathbb{E}_{\mu_K (dx^{B(l)} | y_l)} \left[  \frac{1}{2} \sum_{l=1}^{M} \sum_{n \neq l} \sum_{\substack{i \in B(l) \\ j \in B(n)}} |M_{ij}| x_i ^2   \right] \right) .
\end{align}
Because~$\mu_K(dx^{B(l)} | y_l)$ is a one-dimensional canonical ensemble for each~$l \in [M]$, a similar argument from the proof of Lemma~\ref{l_subadditive_gce} using Corollary~\ref{l_ce_second_moment} yields
\begin{align}
T_{\eqref{e_difference_cg_hamiltonian}} \lesssim \frac{1}{N} \sum_{l=1}^{M} \left( 1 + y_l ^2 \right) = \frac{1}{K} \left(1 + \|y \|_{L^2(Y)} ^2 \right)
\end{align}
The lower bound of~$T_{\eqref{e_difference_cg_hamiltonian}}$ is similarly deduced. \\
\qed

Next, we bound the coarse-grained Hamiltonian~$\bar{H}$ with quadratic functions.

\begin{lemma} \label{l_bound_cg_hamiltonian}
There exists a positive constant~$C$ such that
\begin{align}
-C + \frac{1}{C} \| y \|_{L^2 (Y)}^2 \leq \bar{H}(y) \leq C \left(1 + \|y\|_{L^2 (Y)}^2 \right). 
\end{align}
\end{lemma}

\medskip

\noindent \emph{Proof of Lemma~\ref{l_bound_cg_hamiltonian}.} \ By Lemma~\ref{l_close_cg_hamiltonians}, it suffices to prove
\begin{align} \label{e_aim_bound_cg_aux}
    -C + \frac{1}{C} \| y \|_{L^2 (Y)}^2 \leq \bar{H}_{\text{aux}}(y) \leq C \left(1 + \|y\|_{L^2 (Y)}^2 \right). 
\end{align}
First of all, an application of Lemma~\ref{l_convergence_cg_ham} yields that
\begin{align}
\bar{H}_{\text{aux}}( \mathbf{0}  ) \overset{\eqref{d_aux_cg_ham}}{=} \frac{1}{M} \sum_{l=1}^{M} \bar{H}_K (0) = \bar{H}_K (0) \to \varphi(0) \qquad \text{as } \nu \to \infty.
\end{align}
Thus~$\bar{H}_{\text{aux}}(\textbf{0})$ is uniformly bounded. \\

Next, Lemma~\ref{l_gradient_cg_hamiltonian} implies
\begin{align}
\left| \frac{\partial}{\partial y_l } \bar{H}_{\text{aux}} (\mathbf{0}) \right| &= \left| \frac{1}{N}\mathbb{E}_{\mu_{\text{aux}}(dx|0)}\left[ \sum_{i, j \in B(l)} M_{ij}X_i + \sum_{i \in B(l)} \delta \psi ' (X_i)\right] \right| \\
& \overset{Lemma~\ref{l_ce_moment_estimate}}{\lesssim} \frac{1}{N} \left( 2KR^2 + K \right)  \sim \frac{1}{M}.
\end{align}
In particular the partial derivatives of~$\bar{H}_{\text{aux}}( \textbf{0} )$ are also bounded. \\

Because~$\bar{H}_{\text{aux}}$ is uniformly strictly convex (Theorem~\ref{p_strict_convexity_cg_hamiltonian}), an application of Taylor's theorem establishes the desired inequalities~\eqref{e_aim_bound_cg_aux}.  \\
\qed

\begin{corollary} \label{l_convexity_phi}
The one-dimensional coarse-grained Hamiltonian~$\bar{H}_K$ and its limit~$\varphi$ are strictly convex. In particular, there exists a   constant~$C>0$ such that
\begin{align}
    -C + \frac{1}{C} m^2 \leq \bar{H}_K (m), \varphi(m) \leq C( 1+ m^2).
\end{align}
\end{corollary}

\medskip

The last two ingredients for deducing Proposition~\ref{p_micro_meso} are the following statements. Since the proofs of   Lemma~\ref{l_kappa} and Lemma~\ref{l_92} in~\cite{GrOtViWe09} work without any changes, we just present statements without a proof. \\

\begin{lemma} \label{l_kappa}
Define~$\kappa$ by
\begin{align}
\kappa : = \max \left \{ \langle \Hess H(x) \cdot u, v \rangle : \ u \in \textup{Ran}(NP^* P), v \in \textup{Ran}(\Id_X - NP^* P), |u|=|v|=1    \right \}.
\end{align}
Then, we have $\kappa < \infty$.
\end{lemma}

\medskip

\begin{lemma} [(92), (93) in~\cite{GrOtViWe09}] \label{l_92}
There exists a   constant~$C>0$ such that
\begin{align} 
    \frac{1}{C} \langle \bar{x}, \bar{x} \rangle_{H^{-1}} \leq \frac{1}{N} \langle x , A^{-1} x \rangle_X \leq C \langle \bar{x}, \bar{x} \rangle_{H^{-1}}. \label{e_92}
\end{align}
If~$\bar{x}$ is bounded in~$L^2$, then
\begin{align}
    \left| \langle \bar{x}, \bar{x} \rangle_{H^{-1}} - \frac{1}{N} \langle x , A^{-1} x \rangle_X  \right| \leq \frac{C}{N}. \label{e_93}
\end{align}
\end{lemma}

\subsubsection{Proof of Proposition~\ref{p_micro_meso}} \label{s_proof_micro_meso_conclusion}
In this  section, we prove Proposition~\ref{p_micro_meso}
with the help of Corollary~\ref{a_propagation}. \\

\noindent \emph{Proof of Proposition~\ref{p_micro_meso}.} \ Let us begin with verifying the assumptions of Theorem~\ref{a_two_scale_hydrodynamic}.

\begin{itemize}
    \item (\romannumeral 1) is the same as Lemma~\ref{l_kappa}.
    \item (\romannumeral 2) is a consequence of Theorem~\ref{p_generalized_lsi}.
    \item (\romannumeral 3) is a consequence of Theorem~\ref{p_strict_convexity_cg_hamiltonian}.
    \item (\romannumeral 4) follows from Lemma~\ref{l_ce_second_moment}.
    \item (\romannumeral 5) is a consequence of Lemma~\ref{l_bound_cg_hamiltonian}.
    \item (\romannumeral 6) is the same as~\eqref{e_92} in Lemma~\ref{l_92}
\end{itemize}

The first assumption of (\romannumeral 7) is the same as~\eqref{e_entropy_f0}. To verify the second condition, let us recall~\eqref{d_eta0}. Because~$\zeta_0 \in L^2 (\mathbb{T}^1)$, there is a positive constant~$C$ such that
    \begin{align}
        \| \bar{\eta}_0 ^{\nu} \|_{L^2} \leq C.
    \end{align}
Therefore Lemma~\ref{l_bound_cg_hamiltonian} implies the second condition of (\romannumeral 7) as follows:
\begin{align}
    \bar{H}( \eta_0 ^{\nu}) \leq C (1 + \| \eta_0 ^{\nu} \|_{L^2 (Y)} ^2 ) = C(1 + \| \bar{\eta}_0 ^{\nu} \|_{L^2} ^2 )  \leq C.
\end{align}
\medskip

Let us turn to the proof of~\eqref{e_cor_assumption} in Corollary~\ref{a_propagation}. It holds that
\begin{align} 
0 &\overset{\eqref{e_92}}{\leq} \frac{1}{N} \int (x- N P^t \eta_{0}^{\nu}) \cdot A ^{-1} ( x-N P^t \eta_{0}^{\nu}) f_{0}(x) \mu (dx) \\
& \overset{\eqref{e_92}}{\leq} \int C \| \bar{x} - \bar{\eta}_0 ^{\nu} \|_{H^{-1}}^{2} f_0 (x) \mu(dx) \\
& \ \leq \int 2 C \left( \| \bar{x} - \zeta_0 \|_{H^{-1}}^{2}  + \|   \zeta_0 - \bar{\eta}_0 ^{\nu} \|_{H^{-1}}^{2}  \right) f_0 (x) \mu(dx) \overset{\eqref{e_entropy_f0},\eqref{d_eta0}}{\longrightarrow} 0 \qquad \text{as } \nu \to \infty.
\end{align}

Therefore all assumptions of Corollary~\ref{a_propagation} are verified and we obtain
\begin{align}
\lim_{v \to \infty} \sup_{0\leq t \leq T} \frac{1}{N} \int ( x- N P ^{t} \eta) \cdot A^{-1}( x- N P ^{t} \eta) f(t,x)\mu(dx) = 0.
\end{align}
In particular,~\eqref{e_92} implies
\begin{align}
\lim_{\nu \to \infty} \sup_{0 \leq t \leq T} \int  \| \bar{x} - \bar{\eta} ^{\nu}(t, \cdot) \|_{H^{-1}} f(t, x) \mu_{N, m}(dx) = 0.
\end{align}
\qed

\subsection{Proof of Proposition~\ref{p_meso_macro}} \label{s_proof_meso_macro}

\subsubsection{Auxiliary Lemmas} \label{s_meso_macro_auxiliary}

In this Section, we provide auxiliary statements that will be needed in the proof of Proposition~\ref{p_meso_macro}. The statements are extensions of~\cite{GrOtViWe09}, where the ce without interactions are considered, to the general ce with strong interactions.

\begin{lemma}[Analogue of Lemma 34 in~\cite{GrOtViWe09}] \label{l_lem34} 
 There is a constant~$C>0$ such that
\begin{align}
&\sup_{0 \leq t \leq T} \langle \eta ^{\nu} (t), \eta^{\nu} (t) \rangle_Y \leq C, \label{e_etav_norm_bdd} \\
&\int_0 ^T \left \langle \frac{d \eta^{\nu}}{dt}(t)  , ( \bar{A} )^{-1} \frac{d \eta^{\nu}}{dt}(t)  \right \rangle_Y dt  \leq C.
\end{align}
\end{lemma}
\medskip 

In particular,~\eqref{e_etav_norm_bdd} implies, up to a subsequence, the associated step functions~$\bar{\eta}^{\nu}$ converges  to~$\eta_*$ weak-$*$  in~$L^{\infty}(L^2) = (L^1 (L^2))^*$. Next lemma provides some properties of the function $\eta_*$.

\begin{lemma}[Analogue of Lemma 35 in~\cite{GrOtViWe09}] \label{l_lem35}
Let~$\{ \eta^{\nu} \}_{v=1}^{\infty}$ as in Lemma~\ref{l_lem34} such that it weak* converges to~$\eta_* $ in~$L^{\infty}(L^2) = (L^1 (L^2))^*$. Then~$\eta_*$ satisfies
\begin{align}
\eta_* \in L_t^{\infty}(L_{\theta}^2) , \qquad \frac{\partial \eta_*}{\partial t} \in L_t ^2 (H_{\theta}^{-1}), \qquad \varphi' (\eta_*) \in L_t^2 (L_{\theta}^2).
\end{align}
\end{lemma}
\medskip

The following lemma provides a integral criteria to ensure a function to be a weak solution to the nonlinear parabolic equation.
\begin{lemma}[Analogue of Lemma 36 in~\cite{GrOtViWe09}] \label{l_lem36}
Assume~$\bar{H}$ is convex. Then~$\eta$ satisfies
\begin{align}
\frac{d \eta^{v}}{dt} = - \bar{A}\nabla_Y \bar{H}(\eta ^v)
\end{align}
if and only if for all~$\xi \in Y$ and~$\beta : [0, T] \to [0, \infty) $ smooth,
\begin{align}
\int_0 ^ T \bar{H}(\eta) \beta (t)dt \leq \int_0^{T} \bar{H}(\eta + \xi)\beta(t) dt - \int_0^T \langle \xi, (\bar{A})^{-1} \eta \rangle_Y \dot{\beta}(t) dt.
\end{align}

Similarly, if~$\varphi$ is convex, then~$\zeta$ is a weak solution of
\begin{align}
\frac{\partial \zeta}{\partial t} = \frac{\partial^2}{\partial \theta^2} \varphi' (\zeta) 
\end{align}
if and only if for all~$\xi \in L^2 (\mathbb{T}^1)$ and~$\beta : [0, T] \to [0, \infty)$ smooth,
\begin{align}
\int_0 ^T \int_{\mathbb{T}^1} \varphi( \zeta(t, \theta) )\beta(t) d \theta dt \leq \int_0 ^T \int_{\mathbb{T}^1} \varphi(\zeta(t,\theta)+ \xi(\theta)) \beta(t) d\theta dt - \int_0 ^T \langle \xi(\cdot), \zeta(t, \cdot) \rangle_{H^{-1}} \dot{\beta}(t)dt.
\end{align}
\end{lemma}
\medskip

In addition, we have a uniqueness of the weak solution to the nonlinear parabolic equation~\eqref{e_heat_eqn}.
\begin{lemma}[Lemma 38 in~\cite{GrOtViWe09}] \label{l_lem38}
The weak solution of~\eqref{e_heat_eqn} is unique, if it exists.
\end{lemma}

We shall not provide the proof of Lemma~\ref{l_lem34}, Lemma~\ref{l_lem35}, Lemma~\ref{l_lem36} and Lemma~\ref{l_lem38} as there are only cosmetic differences to that of~\cite[Lemma 34]{GrOtViWe09},~\cite[Lemma 35]{GrOtViWe09},~\cite[Lemma 36]{GrOtViWe09} and~\cite[Lemma 38]{GrOtViWe09}, respectively. The crucial step to relate solutions to the mesoscopic and macroscopic equations is the following lemma. \\

\begin{lemma} [Lemma 37 in~\cite{GrOtViWe09}] \label{l_lem37}
Let~$\{\eta ^{\nu} \}_{\nu=1}^{\infty}$ and~$\eta_*$ as in Lemma~\ref{l_lem35}. Define~$\xi^{\nu}: = \pi (\xi + \eta_* ) - \eta^{\nu}$, where~$\pi$ is the~$L^2$ projection onto~$Y$. Then it holds that
\begin{align}
&\lim_{ \nu \to \infty} \int_0 ^T \bar{H}( \eta^{\nu}(t)) \beta(t)dt \geq \int_0 ^T \int_{\mathbb{T}^1} \varphi( \eta_* (t, \theta) ) \beta(t) d\theta dt, \label{e_lem37_1} \\
& \lim_{\nu \to \infty} \int_0 ^T \bar{H} ( \eta ^{\nu} (t) + \xi^{\nu} (t))\beta(t)dt = \int_0 ^T \int_{\mathbb{T}^1} \varphi( \eta_* (t, \theta)+ \xi(\theta) ) \beta(t) d\theta dt,\label{e_lem37_2} \\
& \lim_{\nu \to \infty} \int_0 ^T \left\langle \xi^{\nu} (t) , (\bar{A})^{-1} \eta^{\nu} (t) \right\rangle_Y \dot{\beta}(t)dt  = \int_0 ^T \langle \xi(\theta), \eta_* (t, \theta) \rangle_{H^{-1}} \dot{\beta}(t) dt. \label{e_lem37_3}
\end{align}
\end{lemma}

\medskip

Since the proof of a statement (108) in \cite{GrOtViWe09} can be adapted to prove~\eqref{e_lem37_3} without any changes, we 
 present the proof of~\eqref{e_lem37_1} and~\eqref{e_lem37_2}. Compared to~\cite{GrOtViWe09} the main difficulty one encounters when deducing~\eqref{e_lem37_1} and~\eqref{e_lem37_2} is the lack of uniform convergence of the coarse grained Hamiltonian~$\bar H_K$ towards~$\varphi$. The key ingredient to solve this problem is Proposition  \ref{l_convergence_cg_ham}, which gives a quantitative convergence of $\bar{H}$ towards $\varphi$.  \\

\noindent \emph{Proof of~\eqref{e_lem37_1} in Lemma~\ref{l_lem37}.} \  Let us write
\begin{align}
\int_0 ^T \bar{H}( \eta^{\nu}(t)) \beta(t)dt  &= \int_0 ^T \left( \bar{H}( \eta^{\nu}(t)) - \bar{H}_{\text{aux}}( \eta^{\nu}(t)) \right) \beta(t)dt \label{e_37_aux_remainder}  \\
& \quad + \int_0 ^T \bar{H}_{\text{aux}}( \eta^{\nu}(t)) \beta(t)dt . \label{e_37_aux_approx}
\end{align}
To begin with, an application of Lemma~\ref{l_close_cg_hamiltonians} followed by~\eqref{e_etav_norm_bdd} yields
\begin{align} \label{e_37_11}
\left| T_{\eqref{e_37_aux_remainder}} \right| \leq \frac{C}{K} \int_0 ^T (1 + \| \eta^{\nu} \|_{L^2} ^2 ) \beta(t)dt \leq \frac{C}{K} \int_0 ^T \beta(t) dt. \label{e_4}
\end{align}
Next, we have
\begin{align}
T_{\eqref{e_37_aux_approx}} &\overset{\eqref{d_aux_cg_ham}}{=} \int_0 ^T  \int_{\mathbb{T}^1} \bar{H}_K (\bar{\eta}^{\nu})  \beta(t) d\theta dt \\
& = \int_0 ^T  \int_{\mathbb{T}^1} \varphi (\bar{\eta}^{\nu}) \beta(t) d\theta dt + \int_0 ^T  \int_{\mathbb{T}^1} \left(\bar{H}_K (\bar{\eta}^{\nu}) - \varphi (\bar{\eta}^{\nu}) \right)  \beta(t) d\theta dt \label{e_3}
\end{align}
Since~$\varphi$ is convex, the functional $f \mapsto \int \int \varphi(f) d\theta dt$ is weak lower semicontinuous with respect to weak-$*$ $L^\infty (L^2)$ topology. Thus, we have
\begin{align} \label{e_37_21}
\lim_{\nu \to \infty} \int_0 ^T  \int_{\mathbb{T}^1} \varphi (\bar{\eta}^{\nu})  \beta(t) d\theta dt  \geq \int_0 ^T \int_{\mathbb{T}^1} \varphi( \eta_* )  \beta(t) d\theta dt . \label{e_1}
\end{align}
Next, we have by Lemma~\ref{l_convergence_cg_ham} that
\begin{align} \label{e_37_22}
\left| \int_0 ^T  \int_{\mathbb{T}^1} \left(\bar{H}_K (\bar{\eta}^{\nu}) - \varphi (\bar{\eta}^{\nu}) \right)  \beta(t) d\theta dt \right| \leq \frac{C}{K} \int_0 ^T (1 + \| \eta^{\nu} \|_{L^2}^2 )\beta(t) dt \overset{\eqref{e_etav_norm_bdd}}{\leq} \frac{C}{K} \int_0 ^T \beta(t) dt
\end{align}
Therefore, plugging the estimations~\eqref{e_37_11},~\eqref{e_37_21} and~\eqref{e_37_22} into~\eqref{e_37_aux_remainder},~\eqref{e_37_aux_approx} and taking~$\nu \to \infty$ yields
\begin{align}
\lim_{\nu \to \infty} \int_0 ^T \bar{H}( \eta^{\nu}(t)) \beta(t)dt \geq \int_0 ^T \int_{\mathbb{T}^1} \varphi(\eta_*) \beta(t) d\theta dt.
\end{align}

\qed

\medskip

\noindent \emph{Proof of~\eqref{e_lem37_2} in Lemma~\ref{l_lem37}.} \ Let us write
\begin{align}
&\int_0 ^T \bar{H} ( \eta ^{\nu} (t) + \xi^{\nu} (t))\beta(t)dt - \int_0 ^T \int_{\mathbb{T}^1} \varphi( \eta_* (t, \theta)+ \xi(\theta) ) \beta(t) d\theta dt \\
& \qquad \overset{\eqref{d_aux_cg_ham}}{=} \int_0 ^T  \int_{\mathbb{T}^1} \bar{H}_K (\bar{\eta}^{\nu} + \bar{\xi}^{\nu})  \beta(t) d\theta dt  - \int_0 ^T \int_{\mathbb{T}^1} \varphi( \eta_* (t, \theta)+ \xi(\theta) ) \beta(t) d\theta dt \\
& \qquad \ = \int_0 ^T  \int_{\mathbb{T}^1} \left( \bar{H}_K (\bar{\eta}^{\nu} + \bar{\xi}^{\nu}) - \varphi(\bar{\eta}^{\nu} + \bar{\xi}^{\nu}) \right) \beta(t) d\theta dt \label{e_37_1} \\
& \qquad \quad \ + \int_0 ^T  \int_{\mathbb{T}^1} \left( \varphi( \bar{\eta}^{\nu} + \bar{\xi}^{\nu}) - \varphi(\eta_* + \xi) \right) \beta(t) d\theta dt. \label{e_37_2}
\end{align}
Let us begin with the estimation of~\eqref{e_37_1}. Lemma~\ref{l_convergence_cg_ham} implies that
\begin{align}
\left| T_{\eqref{e_37_1}}\right| & \leq \frac{C}{K} \int_0 ^T \int_{\mathbb{T}^1} \left( 1+ \| \bar{\eta}^{\nu} + \bar{\xi}^{\nu} \|_{L^2} ^2 \right) \beta(t) d\theta dt \\
& \leq \frac{C}{K} \int_0 ^T \int_{\mathbb{T}^1} (1 + \| \eta_* + \xi\|^2 ) \beta(t) d\theta dt \\
& \leq \frac{C}{K} \int_0 ^T \beta(t) dt \label{e_37_1_est}
\end{align}
Let us turn to the estimation of~\eqref{e_37_2}. Recalling that~$\xi^{\nu}$ is defined by~$ \eta^{\nu} +\xi^{\nu}  = \pi(\xi + \eta_* )$, we have
\begin{align} \label{e_strong_l2_conv}
    \bar{\eta}^{\nu} + \bar{\xi}^{\nu} \to \eta_* + \xi \qquad \text{in } L^2 \text{ for a.e. } t.
\end{align}
A combination of Corollary~\ref{l_convexity_phi} and~\eqref{e_strong_l2_conv} yields
\begin{align} \label{e_pw}
    \int_{\mathbb{T}^1} \varphi( \bar{\eta}^{\nu} + \bar{\xi}^{\nu} ) d\theta \to \int_{\mathbb{T}^1} \varphi(\eta_* + \xi) d\theta \qquad \text{for a.e. } t.
\end{align}
In addition, we have
\begin{align} \label{e_bound}
\left| \int_{\mathbb{T}^1} \varphi( \bar{\eta}^{\nu} + \bar{\xi}^{\nu} ) d\theta \right| \leq  \int_{\mathbb{T}^1} C\left(1 + | \bar{\eta}^{\nu} + \bar{\xi}^{\nu} |^2 \right)d\theta \leq  \int_{\mathbb{T}^1} C\left(1 + | \eta_* + \xi |^2 \right)d\theta \leq C.
\end{align}
Thus the Dominated Convergence theorem applied with~\eqref{e_pw} and~\eqref{e_bound} gives
\begin{align} \label{e_37_last2}
    \lim_{\nu \to \infty} T_{\eqref{e_37_2}} = \int_0 ^T \int_{\mathbb{T}^1} \varphi(\eta_* + \xi) \beta(t) d\theta dt.
\end{align}
Now letting~$\nu \to \infty$ in~\eqref{e_37_1_est} and plugging this with~\eqref{e_37_last2} into~\eqref{e_37_1},~\eqref{e_37_2} yields
\begin{align}
\lim_{\nu \to \infty}  \int_0 ^T \bar{H} ( \eta ^{\nu} (t) + \xi^{\nu} (t))\beta(t)dt = \int_0 ^T \int_{\mathbb{T}^1} \varphi( \eta_* (t, \theta)+ \xi(\theta) ) \beta(t) d\theta dt.
\end{align}
\qed

\medskip

\subsubsection{Proof of Proposition~\ref{p_meso_macro}} In this Section, we prove Proposition~\ref{p_meso_macro}. \\

\noindent \emph{Proof of Proposition~\ref{p_meso_macro}.} \ Because~$\zeta_0 \in L^2(\mathbb{T}^1)$ and~$\bar{\eta}_0^{v} $ converges to~$\zeta_0$ in~$L^2$, we know~$L^2$ norm of~$\bar{\eta}_0 ^v$ is uniformly bounded. Therefore an application of Lemma~\ref{l_lem34} yields that up to a subsequence,
\begin{align}
\bar{\eta}^v \rightharpoonup \eta_* \qquad \text{ weak}* \text{ in } L_t ^{\infty}(L_{\theta}^2) = (L_t ^1 (L_{\theta}^2))^*.
\end{align}
Lemma~\ref{l_lem35} implies the weak* limit~$\eta_*$ satisfies
\begin{align}
\eta_* \in L_t^{\infty}(L_{\theta}^2) , \qquad \frac{\partial \eta_*}{\partial t} \in L_t ^2 (H_{\theta}^{-1}), \qquad \varphi' (\eta_*) \in L_t^2 (L_{\theta}^2).
\end{align}
Next, we have by Lemma~\ref{l_lem36} that
\begin{align} \label{e_result_lem36}
\int_0 ^ T \bar{H}(\eta^v) \beta (t)dt \leq \int_0^{T} \bar{H}(\eta^v + \xi^v)\beta(t) dt - \int_0^T \langle \xi^v, (\bar{A})^{-1} \eta^v \rangle_Y \dot{\beta}(t) dt,
\end{align}
where~$\xi^v := \pi (\xi + \eta_* ) - \eta^v $. By taking the limit in~\eqref{e_result_lem36} and applying Lemma~\ref{l_lem37}, one gets
\begin{align}
\int_0 ^T \int_{\mathbb{T}^1} \varphi( \zeta(t, \theta) )\beta(t) d \theta dt \leq \int_0 ^T \int_{\mathbb{T}^1} \varphi(\zeta(t,\theta)+ \xi(\theta)) \beta(t) d\theta dt - \int_0 ^T \langle \xi(\cdot), \zeta(t, \cdot) \rangle_{H^{-1}} \dot{\beta}(t)dt.
\end{align}
Therefore Lemma~\ref{l_lem36} implies that $\eta_*$ is a weak solution of
\begin{align}
\begin{cases}
\frac{\partial \zeta}{\partial t} = \frac{ \partial^2}{\partial \theta^2}\varphi ' (\zeta), \\
\zeta(0, \cdot) = \zeta_0,
\end{cases}
\end{align}
and by uniqueness (Lemma~\ref{l_lem38}), the sequence~$\{ \bar{\eta}^v \}_{v=1}^{\infty}$ converges to~$\eta_*$. \qed

\newpage

\appendix

\section{Basic Properties of the grand canonical ensemble and the canonical ensemble.}

In this section, we provide auxiliary results which were proved in~\cite{KwMe18},~\cite{KwMe19a},~\cite{KwMe19b} and~\cite{KLM19}. Recall that a generalized  gce~$\mu^{\sigma}_N$ is defined by
\begin{align}
\mu^{\sigma}_N (dx) : = \frac{1}{Z} \exp\left( \sigma \sum_{i=1}^N x_i - H(x) \right) dx.
\end{align}
With the term~$s = (s_i )_i \in \mathbb{R}^N$,~$\sigma \in \mathbb{R}$ models the interaction of the system with boundary values.

\begin{remark}
Again, the ce~$\mu_{N,m}$ can be thought as a conditional probability distribution which emerges from the gce~$\mu_{N}^{\sigma}$ conditioned on the mean spin
\begin{align}
    \frac{1}{N} \sum_{i=1}^{N} x_i = m.
\end{align}
More precisely, we have
\begin{align}
\mu_N^{\sigma}  \left(dx   \mid \frac{1}{N}\sum_{i =1}^{N} x_i =m \right) & = \frac{1}{Z} \mathds{1}_{ \left\{ \frac{1}{N} \sum_{i=1}^{N} x_i =m \right\}}(x) \exp\left(\sigma m N - H(x) \right) \mathcal{L}^{N-1}(dx)  \\
& = \frac{1}{\widetilde{Z}} \mathds{1}_{ \left\{ \frac{1}{N} \sum_{i =1}^{N} x_i =m \right\}}\left(x\right) \exp\left( - H(x) \right) \mathcal{L}^{N-1}(dx) \\
& = \mu_{N,m} (dx).
\end{align}
\end{remark}

The following statement tells that the variance of the mean spin of the modified gce~$\mu^{\sigma}_N$ is well behaved.

\begin{lemma} [Lemma 1 in~\cite{KwMe19b}] \label{l_variance_estimate}
There exists a constant~$C \in (0, \infty)$, uniform in $N$, $s$, and~$\sigma$ such that
\begin{align}\label{e_variance_estimate_mgce}
\frac{1}{C} \leq \frac{1}{N} \var_{\mu^{ \sigma}_N} \left( \sum_{k=1}^N X_k \right) \leq C.
\end{align}
\end{lemma}

Define the free energy~$A_N : \mathbb{R} \rightarrow \mathbb{R}$ by
\begin{align}
A_N(\sigma) : = \frac{1}{N} \log \int_{\mathbb{R}^N} \exp \left( \sigma \sum_{i=1}^N x_i - H(x) \right) dx.
\end{align}
Then the free energy~$A_N$ is uniformly strictly convex.

\begin{lemma} [Lemma 2 in~\cite{KwMe19b}] \label{l_strict_convex_free_energy} There is a constant~$C \in (0, \infty)$, uniform in $N$, $s$, and~$\sigma$ such that
\begin{align}
\frac{1}{C} \leq \frac{d^2}{d\sigma^2} A_N(\sigma) \leq C.
\end{align}
\end{lemma}

Now we relate the external field~$\sigma$ of~$\mu_N ^{\sigma}$ and the mean spin~$m$ of~$\mu_{N,m}$ as follows:

\begin{definition} For each~$m \in \mathbb{R}$, we choose~$\sigma = \sigma_N(m) \in \mathbb{R}$ such that
\begin{align}
\frac{d}{d\sigma} A_N( \sigma) = m.
\end{align}
Denoting~$m_i : = \int x_i \mu_N ^{\sigma} (dx)$ for each~$i \in [N]$, we equivalently get
\begin{align}
m = \frac{d}{d\sigma} A_N ( \sigma) = \frac{1}{N} \frac{ \int_{\mathbb{R}^N} \sum_{i=1}^{N} x_i \exp \left( \sigma \sum_{i=1}^{N} x_i - H(x) \right)dx }{  \int_{\mathbb{R}^N}  \exp \left( \sigma \sum_{i=1}^{N} x_i - H(x) \right) dx }  = \frac{1}{N}\sum_{i=1}^N m_i .
\end{align}
\end{definition}

Let us now introduce the definition of local, intensive and extensive functions.

\begin{definition} [Local, intensive, and extensive functions/ observables]
For a function~$f : \mathbb{R}^{\mathbb{Z}} \to \mathbb{C}$, denote~$\supp f $ by the minimal subset of~$\mathbb{Z}$ with $f(x) = f\left(x^{\supp f} \right)$. We call~$f$ a local function if it has a finite support independent of~$N$. A function~$f$ is called intensive if there is a positive constant~$\varepsilon$ such that $|\supp f | \lesssim N^{1-\varepsilon}$. A function~$f$ is called extensive if it is not intensive.
\end{definition}

For one-dimensional lattice systems the correlations of the gce decay exponentially fast.

\begin{proposition}[Lemma 6 in~\cite{KwMe18}] \label{p_decay_of_correlations_gce}
Let~$f, g : \mathbb{R}^{N} \to \mathbb{R}$ be intensive functions. Then 
\begin{align}
&\left| \cov_{\mu_N ^{\sigma}} \left( f, g \right) \right| \lesssim \|\nabla f\|_{L^2 (\mu_N ^{\sigma})}\|\nabla g\|_{L^2 (\mu_N ^{ \sigma})} \exp \left( -C \text{dist} \left( \supp f , \supp g \right) \right). \label{e_decay_of_correlation_ce}
\end{align}
\end{proposition}

The following moment estimate is a consequence of Proposition~\ref{p_decay_of_correlations_gce}.

\begin{lemma}[Lemma 3.2 in~\cite{KwMe19a}] \label{l_moment_estimate} For each~$k \geq 1$, there is a constant~$C = C(k)$ such that for any smooth function~$f : \mathbb{R}^{\Lambda} \to \mathbb{R}$
\begin{align} \label{e_function_moment_estimate}
\mathbb{E}_{\mu_N^{ \sigma}} \left[ \left|f(X) - \mathbb{E}_{\mu_N^{ \sigma}} \left[f(X) \right] \right|^k \right] \leq C(k) \| \nabla f \|_{\infty}^k.
\end{align}
\end{lemma}

For the ce, the correlations decay exponentially fast with a volume correction term.

\begin{proposition}[Theorem 2.10 in~\cite{KLM19}] \label{p_decay_of_correlations_ce}
Let~$f, g : \mathbb{R}^{N} \to \mathbb{R}$ be intensive functions. There exist constants~$C \in (0, \infty)$ and~$N_0 \in \mathbb{N}$ independent of the external field~$s$ and the mean spin~$m$ such that for all~$N \geq N_0$, it holds that
\begin{align} \label{e_decay_correlations_ce}
\left|\cov_{\mu_{N,m} } \left( f, g \right) \right| \leq C \ \| \nabla f \|_{L^{\infty}(\mu_N^{\sigma})}\| \nabla g \|_{L^{\infty}(\mu_N^{\sigma})}  \left( \frac{ |\supp f | + |\supp g | }{N} + \exp\left(-C\text{dist}\left( \supp f, \supp g \right) \right) \right).
\end{align}
\end{proposition}

Lastly, we introduce the equivalence of observables.

\begin{proposition} [Theorem 2.7 in~\cite{KLM19}] \label{p_equivalence_observables}
Let~$f : \mathbb{R}^N \to \mathbb{R}$ be an intensive function. There are constants~$C \in (0, \infty)$ and~$N_0 \in \mathbb{N}$ independent of the external field~$s$ and the mean spin~$m$ such that for all~$N \geq N_0$, it holds that
\begin{align}
\left| \mathbb{E}_{\mu_N^{\sigma}} \left[ f \right] - \mathbb{E}_{\mu_{N,m}} \left[f \right] \right| \leq C \frac{|\supp f|}{N} \| \nabla f \|_{\infty}.
\end{align}
\end{proposition}

\section{Derivatives of coarse-grained Hamiltonian}
\begin{lemma}[Lemma 1 in~\cite{Me11}] \label{l_change_of_var}
For~$z \in \{w : Pw=0\}$ and~$y \in Y$, let~$H_{(M_{ij})}(z,y)$ be
\begin{align}
H_{(M_{ij})} (z,y) : = \frac{1}{2}\langle z, (Id+(M_{ij}))z \rangle  + \langle z, (M_{ij})NP^* y\rangle + \langle z, s \rangle + \sum_{i=1}^N \delta \psi (z_i + ( NP^* y )_i ), \label{e_def_h_m}
\end{align}
where~$(M_{ij})$ is the interaction matrix (cf.~\eqref{e_d_hamiltonian}). Then
\begin{align}
\bar{H}(y) = \frac{1}{2} \langle y, (Id + P(M_{ij}) NP^* y \rangle_Y + \langle P^* y , s \rangle - \frac{1}{N}\log \int_{Px=0} \exp \left( -H_{(M_{ij})} (x,y) \right) \mathcal{L}(dx).
\end{align}
\end{lemma}
\medskip

\noindent \emph{Proof of Lemma~\ref{l_change_of_var}.} \ For~$x \in \{w : Pw=y\}$ and~$y \in Y$, let~$z = x - NP^*y$. Recalling the identity~$PNP^* = Id_Y$, we have~$z  \in \{w : Pw = 0\}$. We then write the Hamiltonian~$H$ as
\begin{align}
H(x) &= \sum_{i =1 }^{N} \left( \psi (x_i) + s_i x_i +\frac{1}{2}\sum_{j : \ 1 \leq |j-i| \leq R } M_{ij}x_i x_j \right)\\
&= \frac{1}{2} \langle x , (Id + (M_{ij})x \rangle + \langle x , s \rangle + \sum_{i=1}^{N} \psi_i ( x_i) \\
& = \frac{1}{2} \langle z + NP^* y , (Id + (M_{ij}) ) ( z + NP^* y)  \rangle + \langle z + NP^* y , s \rangle + \sum_{i=1}^{N} \psi_i (z_i + (NP^* y)_i ) \\
& = \frac{1}{2} \langle z , (Id + (M_{ij})) z \rangle + \langle z, NP^* y \rangle + \langle z , (M_{ij}) NP^* y \rangle + \frac{1}{2} \langle NP^* y, (Id + (M_{ij}) ) NP^* y \rangle \\
& \quad + \langle z, s \rangle + \langle NP^* y, s \rangle + \sum_{i=1}^{N} \psi_i (z_i + (NP^* y)_i ). \label{e_hamiltonian_translation}
\end{align}
Because~$Pz = 0$, we have
\begin{align} \label{e_hamiltonian_cancellation}
\langle z , NP^* y \rangle = N \langle Pz, y \rangle_Y = 0.
\end{align}
It also holds by~$PNP^* = Id_Y$ that
\begin{align} \label{e_hamiltonian_reformuation}
\frac{1}{2} \langle NP^* y , (Id +(M_{ij}))NP^* y \rangle &= \frac{N}{2} \langle y , PNP^* y + P(M_{ij}) NP^ * y \rangle_Y  \\
& = \frac{N}{2} \langle y , (Id + P(M_{ij})NP^* y \rangle_Y.
\end{align}
Plugging~\eqref{e_hamiltonian_cancellation} and~\eqref{e_hamiltonian_reformuation} into~\eqref{e_hamiltonian_translation} yields
\begin{align}
H(x) =  N \left( \frac{1}{2} \langle y, (Id + P(M_{ij}) NP^* y \rangle_Y + \langle P^* y , s \rangle \right) +  H_{M_{ij}} (z, y) , 
\end{align} 
and hence
\begin{align}
\bar{H}(y) &= - \frac{1}{N} \log \int_{Px = y} \exp \left(- H(x) \right)\mathcal{L} (dx) \\
&=  \frac{1}{2} \langle y, (Id + P(M_{ij}) NP^* y \rangle_Y + \langle P^* y , s \rangle - \frac{1}{N}\log \int_{Pz=0} \exp \left( -H_{(M_{ij})} (z,y) \right) \mathcal{L}(dz).
\end{align}

\qed

Next, we compute the derivatives of the coarse-grained Hamiltonian~$\bar{H}$ using Lemma~\ref{l_change_of_var}.

\begin{lemma} \label{l_gradient_cg_hamiltonian}
For each~$l \in \{1, \cdots, M\} $, it holds that
\begin{align}
\frac{\partial}{\partial y_l } \bar{H} (y) = \frac{1}{M} y_l + \frac{1}{N} \sum_{j \in B(l)} s_j +  \frac{1}{N}\mathbb{E}_{\mu_{N,m}(dx|y)}\left[ \sum_{i=1}^N \sum_{j \in B(l)} M_{ij}X_i + \sum_{i \in B(l)} \delta \psi ' (X_i)\right].
\end{align}
\end{lemma}

\medskip 

\noindent \emph{Proof of Lemma~\ref{l_gradient_cg_hamiltonian}.} \ Recall the inner product~$\langle \cdot, \cdot \rangle_Y$ is given by
\begin{align}
\langle x, y \rangle_Y : = \frac{1}{M} \sum_{l =1}^M x_l y_l.
\end{align}
First of all, noting that
\begin{align}
\langle y, P(M_{ij})NP^* y \rangle_Y =\frac{1}{N} \sum_{l, n = 1}^{M}\sum_{i \in B(l), j \in B(n)} M_{ij} y_l y_n,
\end{align}
we have
\begin{align}
\frac{\partial}{\partial y_l} \left( \frac{1}{2} \langle y , (Id+ P(M_{ij})NP^* )y \rangle_Y \right) = \frac{1}{M} y_l + \frac{1}{N} \sum_{n=1}^M \sum_{i \in B(l), j \in B(n)} M_{ij} y_n. \label{e_gradient_1st}
\end{align}
Second, a direct calculation yields
\begin{align}
\frac{\partial}{\partial y_l} \langle P^* y , s \rangle = \frac{\partial}{\partial y_l} \left( \frac{1}{N} \sum_{k =1}^{M} \sum_{j \in B(k)} s_j y_k  \right) = \frac{1}{N} \sum_{j \in B(l)} s_j. \label{e_gradient_2nd}
\end{align}
Lastly, differentiating~\eqref{e_def_h_m} yields
\begin{align}
&\frac{\partial}{\partial y_l}\left( H_{(M_{ij})} (x,y) \right) \\
&\qquad  = \frac{\partial}{\partial y_l} \left( \langle x, (M_{ij})NP^* y \rangle \right)  + \frac{\partial}{\partial y_l} \left(\sum_{i=1}^N \delta \psi ( x_i + (NP^* y)_i ) \right) \\
&\qquad  = \sum_{i=1}^N \sum_{j \in B(l)} M_{ij} x_i + \sum_{i \in B(l)} \delta \psi ' (x_i + (NP^* y)_i ) \\
&\qquad  = \sum_{i=1}^N \sum_{j \in B(l)} M_{ij} (x_i + (NP^* y) _i ) - \sum_{i=1}^N \sum_{j \in B(l)} M_{ij}  (NP^* y) _i + \sum_{i \in B(l)} \delta \psi ' (x_i + (NP^* y)_i ) \\
&\qquad  = \sum_{i=1}^N \sum_{j \in B(l)} M_{ij} (x_i + (NP^* y) _i )  + \sum_{i \in B(l)} \delta \psi ' (x_i + (NP^* y)_i ) - \sum_{n=1}^M \sum_{ i \in B(n), j \in B(l)}M_{ij} y_n.
\end{align}
As a consequence we obtain
\begin{align}
&\frac{\partial}{\partial y_l} \left( - \frac{1}{N}\log \int_{Px=0} \exp \left( -H_{(M_{ij})} (x,y) \right) \mathcal{L}(dx) \right) \\
& \qquad = \frac{1}{N} \frac{\int_{Px=0} \frac{\partial}{\partial y_l} \left( H_{(M_{ij})} (x,y) \right) \exp \left( -H_{(M_{ij})} (x,y) \right) \mathcal{L}(dx)}{\int_{Px=0} \exp \left( -H_{(M_{ij})} (x,y) \right) \mathcal{L}(dx)} \\
& \qquad = \frac{1}{N} \mathbb{E}_{\mu_{N,m}(dx|y)}\left[ \sum_{i=1}^N \sum_{j \in B(l)} M_{ij}X_i + \sum_{i \in B(l)} \delta \psi ' (X_i) - \frac{1}{N} \sum_{n=1}^M \sum_{i \in B(n), j \in B(l)} M_{ij} y_n \right]. \label{e_gradient_3rd}
\end{align}

Combining~\eqref{e_gradient_1st},~\eqref{e_gradient_2nd} and~\eqref{e_gradient_3rd} with symmetry of~$M_{ij}$, i.e.,~$M_{ij} = M_{ji}$, we have the desired equation
\begin{align}
\frac{\partial}{\partial y_l } \bar{H} (y) = \frac{1}{M} y_l + \frac{1}{N} \sum_{j \in B(l)} s_j +  \frac{1}{N}\mathbb{E}_{\mu_{N,m}(dx|y)}\left[ \sum_{i=1}^N \sum_{j \in B(l)} M_{ij}X_i + \sum_{i \in B(l)} \delta \psi ' (X_i)\right].
\end{align}
\qed

The second derivatives of the coarse-grained Hamiltonian~$\bar{H}$ follow from a similar calculations.

\begin{lemma}[Lemma 2 in~\cite{Me11}] \label{l_hessian_cg_hamiltonian} For~$1 \leq l, n \leq M$, we have
\begin{align}
\left( \Hess _Y \bar{H} (y) \right)_{ln} &= \delta_{ln} + \delta_{ln} \frac{1}{K} \int \sum_{i \in B(l)} \delta \psi_i '' (x_i) \mu_{N, m} (dx | y) + \frac{1}{K} \sum_{i \in B(l), j \in B(n)}M_{ij} \\
& \quad - \frac{1}{K} \cov_{\mu_{N, m}(dx|y)} \left( \sum_{j \in B(l)} \left( \sum_{i=1}^{N} M_{ij} X_i + \delta \psi_j ' (X_j) \right) ,   \right. \\
& \qquad \qquad \qquad \qquad \qquad \qquad  \left. \sum_{j \in B(n)} \left( \sum_{i=1}^{N} M_{ij} X_i + \delta \psi_j ' (X_j) \right)   \right).
\end{align}

\end{lemma}

\medskip

\section{Criteria for the logarithmic Sobolev inequality} \label{a_lsi_criteria}
In this section we state several standard criteria for deducing a LSI. For proofs we refer to the literature. For a general introduction and more comments on the LSI we refer the reader to~\cite{Led01,Led01a,Roy99,BaGeLe14}. 

\begin{theorem}[Tensorization Principle~\cite{Gro75}] \label{a_tensorization}
Let~$\mu_1$ and~$\mu_2$ be probability measures on Euclidean spaces~$X_1$ and~$X_2$ respectively. Suppose that~$\mu_1$ and~$\mu_2$ satisfy LSI$(\rho_1)$ and LSI$(\rho_2)$ respectively. Then the product measure~$\mu_1 \bigotimes \mu_2$ satisfies LSI$(\rho)$, where~$\rho = \min \{\rho_1, \rho_2\}$.
\end{theorem}

\begin{theorem}[Holley-Stroock Perturbation Principle~\cite{HolStr87}] \label{a_holley_stroock}
Let~$\mu_1$ be a probability measure on Euclidean space~$X$ and~$\delta \psi : X \to \mathbb{R}$ be a bounded function. Define a probability measure~$\mu_2$ on~$X$ by
\begin{align}
\mu_2 (dx) : = \frac{1}{Z} \exp \left( - \delta \psi (x) \right) \mu_1 (dx).
\end{align}
Suppose that~$\mu_1$ satisfies LSI$(\rho_1)$. Then~$\mu_2$ also satisfies LSI with constant
\begin{align}
\rho_2 = \rho_1 \exp \left( - \text{osc } \delta \psi \right),
\end{align}
where~$\text{osc } \delta \psi : = \sup \delta \psi - \inf \delta \psi $.
\end{theorem}
\begin{theorem}[Bakry-\'{E}mery criterion~\cite{BaEm85}]\label{a_bakry_emery} Let~$X$ be a~$N$-dimensional Euclidean space and $H \in C^2 (X)$. Define a probability measure~$\mu$ on~$X$ by
\begin{align}
\mu(dx): = \frac{1}{Z} \exp\left( - H(x)\right) dx.
\end{align}
Suppose there is a constant~$\rho>0$ such that~$\Hess H \geq \rho$. More precisely, for all~$u, v \in X$,
\begin{align}
\langle v, \Hess H(u) v \rangle \geq \rho |v|^2.
\end{align}
Then~$\mu$ satisfies LSI$(\rho)$.
\end{theorem}
\begin{theorem}[Otto-Reznikoff Criterion~\cite{OttRez07}] \label{a_otto_reznikoff} Let~$X = X_1 \times \cdots \times X_N$ be a direct product of Euclidean spaces and~$H \in C^2 (X)$. Define a probability measure~$\mu$ on~$X$ by
\begin{align}
\mu(dx) : = \frac{1}{Z} \exp \left( - H(x) \right) dx.
\end{align}
Assume that
\begin{itemize}
\item For each~$i \in \{1, \cdots, N\}$, the conditional measures~$\mu (dx_i | \bar{x}_i )$ satisfy LSI$(\rho_i)$.
\item For each~$ 1 \leq i \neq j \leq N$ there is a constant~$\kappa_{ij} \in (0, \infty)$ with
\begin{align}
\left| \nabla_i \nabla_j H(x) \right| \leq \kappa_{ij}. \qquad \text{for all } x \in X.
\end{align}
Here,~$| \cdot |$ denotes the operator norm of a bilinear form.
\item Define a symmetric matrix~$A= (A_{ij})_{1 \leq i, j \leq N}$ by
\begin{align}
A_{ij} = \begin{cases} \rho_i, \qquad &\text{if } i=j \\
- \kappa_{ij}, \qquad &\text{if } i \neq j
\end{cases}.
\end{align}
Assume that there is a constant~$\rho \in (0, \infty)$ with
\begin{align}
A \geq \rho \Id,
\end{align}
in the sense of quadratic forms.
\end{itemize}
Then~$\mu$ satisfies LSI$(\rho)$.
\end{theorem}

\begin{theorem}[Two-Scale Criterion for LSI~\cite{GrOtViWe09}] \label{a_two_scale}
Let~$X$ and~$Y$ be Euclidean spaces. Consider a probability measure~$\mu$ on~$X$ defined by
\begin{align}
\mu(dx) : = \frac{1}{Z} \exp\left(-H(x)\right)dx.
\end{align}
Let~$P : X \to Y$ be a linear operator such that for some~$N \in \mathbb{N}$,
\begin{align}
PNP^* = \Id_Y.
\end{align}
Define
\begin{align}
\kappa : = \max \left \{ \langle \Hess H(x) \cdot u, v \rangle : \ u \in \text{Ran}(NP^* P), v \in \text{Ran}(\Id_X - NP^* P), |u|=|v|=1    \right \}.
\end{align}
Assume that
\begin{itemize}
\item $\kappa < \infty$ 
\item There is~$\rho_1 \in (0, \infty)$ such that the conditional measure~$\mu(dx | Px=y)$ satisfies LSI$(\rho_1)$ for all~$y \in Y$.
\item There is~$\rho_2 \in (0, \infty)$ such that the marginal measure~$\bar{\mu}= P_{\#}\mu$ satisfies LSI$(\rho_2 N)$.
\end{itemize}
Then~$\mu$ satisfies LSI$(\rho)$, where
\begin{align}
\rho : = \frac{1}{2} \left( \rho_1 +\rho_2 + \frac{ \kappa^2}{\rho_1} - \sqrt{ \left( \rho_1 + \rho_2 + \frac{\kappa ^2 }{\rho_1} \right)^2 - 4 \rho_1 \rho_2 }    \right) >0.
\end{align}
\end{theorem}

\section*{Acknowledgment}
This work is dedicated to Thomas M.~Liggett (March 29, 1944 – May 12, 2020). He was a great friend, mentor, teacher and colleague to us all at UCLA.\\

This research has been partially supported by NSF grant DMS-1407558. The authors want to thank Felix Otto and H.T.~Yau for bringing this problem to their attention. The authors are also thankful to many people discussing the problem and helping to improve the preprint. Among them are Tim Austin, Frank Barthe, Marek Biskup, Pietro Caputo, Jean-Dominique Deuschel, Max Fathi, Andrew Krieger, Michel Ledoux, Sangchul Lee, Felix Otto, Daniel Ueltschi, and Tianqi Wu. The authors want to thank Marek Biskup, UCLA and KFAS for financial support.

\bibliographystyle{alpha}
\bibliography{bib}

\end{document}